\documentclass[a4paper,10pt]{article}
\usepackage{hyperref}

\usepackage[utf8]{inputenc}

\usepackage[english ]{babel}
 
\usepackage{xcolor}
\usepackage{bbding}

\usepackage[utf8]{inputenc}

\usepackage{tikz-cd}

\usepackage{amsfonts}
\usepackage{amssymb}
\usepackage{graphicx}
\usepackage{mathtools}
\usepackage{skak} 
\usepackage{foekfont}

\usepackage[T1]{fontenc}

\usepackage{amsmath,mathdots}

\usepackage[utf8]{inputenc}
\usepackage{devanagari}

\usepackage{mathtools}
\usepackage{mathdots}
\usepackage{tikz}

\definecolor{amber(sae/ece)}{rgb}{1.0, 0.49, 0.0}
\newfont{\rsfsten}{rsfs10 scaled 1200}

\usepackage{MnSymbol}
\usepackage{tikz}
\usepackage{amscd}
\usepackage{MnSymbol}
\usepackage{wasysym}

\usepackage{mathrsfs}

\newcommand*{\rom}[1]{\expandafter\@slowromancap\romannumeral #1@}

\usepackage{xcolor}

\usepackage{wrapfig}
\newcommand{\tightunderset}[2]{%
  \mathop{#2}\limits_{\vbox to .3ex{\kern-0.95ex\hbox{$#1$}\vss}}}

\usepackage[utf8]{inputenc}

\usepackage{CJKutf8}

\usepackage{pmboxdraw}

\usepackage{stmaryrd}
\usepackage{url}
\usepackage{harmony}
\usepackage{ifsym}

\usepackage{eufrak}

\usepackage{rotating}

\usepackage{dingbat}

\title {Scalar Curvature, Injectivity Radius and
Immersions with Small Second Fundamental
Forms}

\author{Misha Gromov}

\begin{document}

\maketitle

\hspace {70mm}{\it To 80th birthday of B.L.}

\begin{abstract}

We  prove in special cases  the following.

 $\bullet_{Sc}$   Bounds on the {\it  injectivity   radii} of "topologically complicated" Riemannian $n$-manifolds $X$, where the scalar curvatures of $X$   are bounded from below,
 $Sc(X)\geq \sigma>0$.  

 $\bullet_{curv}$ Lower bounds on {\it focal radii}  of  smooth  immersions from   $k$-manifolds, e.g. homeomorphic to the $k$-torus,   to  certain  Riemannian manifolds   of dimensions  $n=k+m$, e.g. 
to the cylinders  $S^{n-1} \times  \mathbb R^1$.

$\bullet_{mean}$ Topological lower bounds on the mean curvatures of domains in Riemannian manifolds. e.g. in the Euclidean $n$-space $\mathbb R^n$.

At the present moment, our results are  limited by the  {\it spin condition} and the   {\it $n\leq 8$ restriction.}

\end {abstract}

  \tableofcontents


\section {Definitions,  Conjectures and Illustrative  Examples}\label {1}

{\it  Terminology and Notation.}  
Given a Riemannian
 manifold   $X =(X,g)$, possibly   non-complete and/or
 with a
 boundary let 
 $$inrad_x(X)=dist (x, \partial_{metr}   X), \mbox { } x\in X, $$
 be the radius $R$ of the maximal open ball $B=B_x(R)\subset X)$, which doesn't intersect the boundary $\partial X\subset X$   as well    the "infinity" of $X$, i.e. the boundary of $X$ in the metric completion of $X$, where  our $\partial_{metr}X$ stands for the    topological boundary of the interior of $X$ in the metric completion of $X$. 
 
  For instance $inrad_x(X)=dist(x, \partial X)$ for (metrically)   complete manifolds with boundaries.

  {\it The  conjugacy radius} of   $X$,
$$ conj.rad_x(X)\leq inrad_x(X), 
$$   
is the maximal $r$, such that  the map $\exp_x$   is {\it a local diffeomorphisms} on the open $r$-ball  $BT_x(r)\subset T_x(X),$
 where, in fact, "$\exp_x$" is defined on the ball $BT_x(R=inrad_x(X))\subset T_x(X)$.

 {\it The  injectivity  radius}   
 $$inj.rad_x(X)\leq conj.rad_x(X)$$ 
 is the maximal $r$,  such that the map $\exp_x:T_x(X)\to X$ is a   {\it diffeomorphism}    from the open 
$r$-ball $BT_x(r)\subset  T_x(X)$ to its   image in $X$, where, in fact, "$\exp_x$" is defined  the ball $BT_x(R=inrad_x(X))$.

Observe,   that the {\it injectivity radius of the universal covering} $\tilde X$ of $X$  {\it mediates }between  the two above radii, 
$$conj.rad_x(X)\geq inj.rad_x(\tilde X)\geq inj.rad_x(X),$$
where $ inj.rad_x(\tilde X)$ stands for $inj.rad_{\tilde  x}(\tilde X)$, where
$\tilde x\in\tilde X$ is a lift of $x$ to
$\tilde X$ and  where the latter radius is independent of the lift.

 {\it Remark.} The ratio $conj.rad(X)/inj.rad(\tilde X)$   can be arbitrarily large. For instance, one can construct   
 
 {\sf Riemannian metrics $g_i$, $i=1,2,...,$, on all compact $3$-manifolds $X$, such that  
 $diameter(X,g_i)\leq 1$ and  $sect.curv(g_i)\leq 1/i$, which makes $conj.rad_x (X,g_i)\geq \sqrt i$ for all $x\in X$.}
 
 (Probbaly, similar metrics exist on all manifolds of dimensions $n\geq 3$.)
 \vspace {1mm}

 {\it The focal radius} of a   smooth immersion\footnote{"Immersions" are {\it locally} (smooth) {\it embeddings}, or, equivalently, smooth maps  with  {\it non-vanishing  differentials} (on the non-zero tangent vectors).}  map 
   $f:Z\to X$, denoted 
 $$foc.rad_z(Z)=foc.rad_z(Z\overset{f}\hookrightarrow X) $$  
is the maximal $r$, for which 
 the differential of the  {\it normal exponential  map } 
$$\exp^\perp_z=\exp_{f(z)}^\perp:BT_{f(z)}^\perp(Z)\to X,$$
is  an immersion.
 \vspace {1mm}

 {\it The  scalar curvature} of  a Riemannian manifold $X=(X,g)$, denoted 
 $$Sc=Sc(X,x)=Sc(X,g)=Sc(g)=Sc_g(x),$$
  is a   function on $X$, which is equal to
 {\sl  the sum of the values  of the   sectional curvatures  at the $n(n-1)$  ordered bivectors  of  an orthonormal frame in $X$, 
 $$Sc(X,x)=\sum_{i, j} sect.curv_x(\tau_{i,j}),\mbox {  }  i\neq j=1,...,n,$$}\vspace {1mm}
 where this sum doesn't depend on the choice of the frame by the Pythagorean theorem
 and where a  {\it characteristic  property} of the scalar curvature is  {\it  additivity under Cartesian-Riemannian products:} 
 $$Sc(X_1\times X_2, g_1\plus g_2)=Sc(X_1,g_1)+ Sc(X_2, g_2).$$

{\it Examples.}  (a) The scalar curvature  of the   $n$-sphere of radius $R$  is
$$Sc(S^n(R)) =n(n-1) sect.curv(S^n(R))=\frac {n(n-1)}{R^2}.$$

(b) The scalar curvature of   products of spheres by  Euclidean spaces,  
$$X=S^m(R)\times \mathbb R^k $$
satisfy
$$Sc(X) =\frac {m(m-1)}{R^2},$$
where, observe, the sectional  curvatures of these products are pinched between $0$  and 
$\frac {Sc(X)}{m(m-1)}$,
$$0\leq sect.curv(X)\leq \frac {Sc(X)}{m(m-1)}= \frac {1} {R^2},$$
and their conjugacy and  injectivity radii are 
$$conj.rad (X)=inj.rad(X) =inj.rad(S^m(R))=R\pi.$$

 (c)  If $X\subset \mathbb R^{n+1}$ is a smooth cooriented hypersurface  with principal curvatures  
$\lambda_i=\lambda_i(x)$, $i=1, ..., n=dim(X)$,  then, by {\it the Gauss formula},  the scalar curvature of the 
Riemannian  metric $g$ on $X$ induced from the Euclidean/Riemannian one is 
$$Sc_g(x)=\sum_{i\neq j} \lambda_i\lambda_j=\left (\sum_i\lambda_i\right)^2-\sum_i \lambda_i^2.$$

\vspace{1mm}

   {\it \textbf {Leon Green's Conjugate Radius  Inequality}} (theorem 5.1  in [G1963]). 
 {\sf {The {\it minimal conjugacy radius} of a complete Riemannian  manifold $X$
   with {\it finite volume} and with the Ricci curvature   bounded from below, $Ricci(X)  \geq R>-\infty$\footnote {In the original paper  [G1963], $X$ is assumed compact, but  if $X$ is non-compact and geodesically complete, then $Vol(X)<\infty$ and  $\inf Ricci(X)>-\infty$  suffice for the proof. } of {dimension} $n$   is  bounded 
    by the {\it mean scalar curvature of $X$}  as follows: 
  $$\inf_{x\in X} (conj.rad_x (X)^2\cdot {vol(X)^{-1} \int _XSc(X,x)dx}\leq \pi^2\cdot n(n-1),$$
 where
    the equality holds {\it only} for (necessarily compact)  manifolds with constant sectional curvatures $>0$}  and where, observe, 
$$  vol(x)^{-1} \int_X Sc(X,x)dy\geq inf Sc(X)= \inf_ {x\in X}Sc(X,x).$$
Consequently, 
$${inf Sc(X)\cdot   ({\sf inf}conj.rad (X))^2}\leq n(n+1) \pi^2 \leqno { \color {teal!24!black}  [Sc|conjr]_n}$$}
 {\it for all compact Riemannian  $n$-manifolds $X$.}\vspace {1mm}

Below is a conjectural  topological sharpening  of $[Sc|conjr]_n$.

\vspace {1mm}

 {\color {red!5!black} \textbf {1.A.}} \textbf { \color  {red!65!black} $[Sc|conjr]_{n-k}$-Conjecture.} {\sf Let $X$ be a complete Riemannian $n$-manifold, such that

 {\it there exists a  homology  class  $h\in H_k(X)$, such that  {\it no non-zero multiple} of $h$ is  representable by a continuous map from
  a compact oriented  $k$-manifold  with positive scalar curvature.}\footnote 
{Products $Y^m\times \mathbb T^k$ are basic examples of such  $X$, see [SY1979'], [GL1980]}, }

\hspace {-6mm} where this is   expressed   in writing by  the inequality
$$Sc(X)  \ngtr_{H_{k,\mathbb Q}} 0.$$
(Informally speaking,  
 the { \it  "$\mathbb Q$-homological  scalar curvature"}   of $X$ is   not {\it  "fully positive in dimension $k$}".)

Then 
{\sf $$ {inf Sc(X)\cdot   (conj.rad (X))^2}\leq \pi^2{m(m-1) }, m=n-k,\leqno { \color {teal!24!black}  [Sc|conjr]_{n-k}}$$
Moreover, if $X$ is  {\it connected non-compact complete}, then 
 $$inf Sc(X)\cdot   ({\sf \inf} conj.rad (X))^2\leq \pi^2(m-1)(m-2),\footnote {If $X$ has  {\it subexponential volume growth}  and $\inf Ricci(X)>\infty$, then the  inequality  ${inf Sc(X)\cdot   (conj.rad (X))^2} \leq  \pi^2 n(n-1)$ follows by 
 Green's integration  argument.

 But  no bound $inf Sc(X)\cdot (inj.rad(X))^2  < C =C_n$ (not even, for $C=\infty$)  
  is known for   general  non-compact complete  manifolds,  except for  $n\leq 5$, where quantifying   the  uniform contractibility  argument from [CL2020] (compare   
  with [CL 2020], 3.10.3 in  [Gr 2021] and 1.1.A in he next secion) shows that  $inf Sc(X)\cdot (inj.rad(X))^2  \leq  10^{100}  $  for complete  Riemannian manifolds $X$ of dimensions $n\leq 5$.}$$

{\it Two Obvious  (conjectural) Corollaries.} Since $inj.rad\leq conj.rad$   and $$conjrad(X)\geq{ \pi \over\sqrt { \max (0, {\sf sup}sect.curv(X))}},$$ where   $ {\sf sup}sect.curv(X)$ is the supremun of the sectional curvatures over all tangent planes in  $X$,
the  $[Sc|conjr]_{m}$-inequality
\hspace {0mm}$ {inf Sc(X)\cdot   (conj.rad (X))^2}\leq \pi^2m(m-1)$
 implies that 
  $${inf Sc(X)\cdot   (inj.rad (\tilde X))^2}\leq \pi^2{m(m-1) }  \leqno { \color {teal!24!black}  [ScInjr]_m}$$
where $\tilde X$ is the universal covering of $X$,
and that
$${inf Sc(X)\over    {\sf sup}sect.curv(X)}\leq  m(m-1),\leqno { \color {teal!24!black}  [Sc/sect]_m} $$}

  Although,  both inequalities  remain conjectural, we shall show  in section 4 
  that  some manifolds $X$  must satisfy at last one of the  two. 
  
  Below is an  instance of this.

 \vspace {1mm}

\textbf {1.B. Two  Examples.} {\sf Let $X$ be a  Riemannian manifold of dimension $n=m+k$, let $Z$ be a closed connected orientable manifold $Z$, which  admits    a metric with {\it non-positive  sectional curvature} and 
let 
$$\phi: X \to  Z,$$
be a continuous  map. 

(For instance, $X\underset{homeo} = Y\times \mathbb T^k$, where  $\mathbb T^k$ is the $k$-torus  and $\phi:(y,t)\mapsto t$.)

Let the scalar curvature   of $X$ be bounded from below by 
 $$Sc(X)\geq   \sigma>0.$$
 
 Let the sectional   curvature of $X$ be bounded from above by
$$sect.curv(X)\leq \kappa, \mbox {  }  \kappa>0, \footnote {This is a rather annoying (possibly redundant for complete $X$) condition, where, necessarily $ \kappa\geq {\sigma\over n(n-1)}$.}$$ 

\textbf {1.B(i). Compact Case.} Let $X$ be compact without boundary, let $dim(Z)=k$   and let the $k$-dimensional induced homology  homomorphism 
$$\phi_\ast: H_k(X)\to H_k( Z)=\mathbb Z $$ 
doesn't vanish.

Let $$m(m-1)\kappa \leq \sigma.$$ 

{\it If 

 {\it  $\bullet_{spin}$  \hspace {1mm}  either $m\leq 3$ or   the universal covering  $\tilde X$ of  $X$ is {\it spin}\footnote{This means that  the restrictions of the  tangent bundle $T(\tilde X)$  to  surfaces in $\tilde X$ are trivial, for which       vanishing of the second homotopy group $\pi_2(X)$ is sufficient.}

 $\bullet_{\leq 8}$  \hspace {4mm}    $n=dim(X)\leq 8$,}

then the universal covering of $X$ satisfies
$$ inj.rad(\tilde X)\leq \pi \sqrt {1\over \kappa}.$$ }}

{\it Remark.}  Albeit highly constrained, this inequality is 
non-vacuous.  For instance, when applied to 
$$X_\lambda=(S^m\times Z, g_{S^m} +\lambda^2g_Z), \mbox { }  \lambda\to \infty,$$ 
it shows  that  {\sl compact Riemannian  manifolds  $Z=(Z, g_Z)$ with $Sc(g_Z)>0$   of dimension $k=8-m$
 admit no metrics with $sect.curv\leq 0$}. (Of course, this is known for all $k$.)

\textbf {1.B(ii). Non-Compact Case.} {\sf  Let $dim(Z)=k-1$, let   $\phi : X\to Z$  be a fibration (e.g.  a topological  bundle of balls),  let   $\hat \psi: \tilde  Z\to X$ be a "lift" of the universal covering of $Z$ to $X$,     i.e. 
$\phi\circ \hat\psi =\pi$}. (For instance, if the fibration $\phi: X\to Z$  admit a section $Z\hookrightarrow  X$, then  the universal covering of the image of this section may be taken  for $\hat \psi$.) 
Let the scalar and the sectional curvatures    of $X$ be bounded as earlier by  
 $$Sc(X)\geq   \sigma\mbox { and }  sect.curv(X)\leq \kappa.$$
 
 Let $\hat R$ be the distance from he image of $\hat \psi$  to the metric boundary of $X$, 
 $$\hat R=dist_{metr} (\hat \psi(\tilde Z), \partial X),$$, 
  (where, by definition,  $\hat R=\infty$  for geodesically complete  $X$),  let 
  $$\hat U=U_{\frac {5}{6}\hat R}(\hat \psi(\tilde Z))\subset X$$
   be the $\frac {5}{6}\hat R$-neighbourhood of  $ \hat \psi(\tilde Z)\subset X$, and let
  $$\hat\sigma =\sigma - {36(n-1) \pi^2\over n\hat R^2}\geq  m(m-1) \kappa.$$

{\it If  

 {\it  $\bullet_{spin}$  \hspace {1mm}  either $m\leq 3$ or   the universal covering of  $X$ is {\it spin}, 
 
 $\bullet_{\leq 8}$  \hspace {4mm}    $n=dim(X)\leq 8$,}

then he injectivity  radius of the universal covering $\tilde X$ of $X$ is bounded in $U$ as follows:
$$\inf_{x\in \hat U} inj.rad_x(\tilde X)\leq \pi \sqrt {m(m-1) \over \hat\sigma}.$$ }

{\textbf {1.C. Corollary: Focal Radius Bound.} {\sf Let $X$ be a complete  Riemannian $n$-manifold such that
$$Sc(X)\geq \sigma > 0, \mbox { }
sect.curv(X) \leq \kappa  \mbox  {  and  } inj.rad (X)\geq r\geq  \pi/\sqrt\kappa. $$ 
Let $Z$  be a  compact  $(n-m)$-dimensional manifold, which admits a metric with $sect.curv\leq 0$ and let
$Z\hookrightarrow X$ be a  smooth  immersion. 

If 
{\it If  $m>0$ and   

 {\it  $\bullet_{spin}$  \hspace {1mm}  either $m\leq 4$ or   the universal covering of  $X$ is {\it spin}, 
 
 $\bullet_{\leq 8}$  \hspace {4mm}    $n=dim(X)\leq 8$,} 
 
 Then
$$foc.rad  (Z)\leq 6\pi\sqrt{n-1 \over n}\sqrt {1\over  \sigma -(m-1)(m-2)\kappa}$$}}}

{\it  Proof.}  Let $T^\perp(Z) \to Z$ be the  normal  bundle of $Z\hookrightarrow X$ and let $BT^\perp(R) \subset  T^\perp(Z) $ be the open $R$-ball subbundle. 

 Observe that the exponential map  $\exp^\perp :BT^\perp(R) \to X$ is defined  for 
 $R\leq foc.rad(Z)$ and is an immersion for these $R$.
 
 Then  1.B(ii) applies to $X^\perp=BT^\perp(foc.rad(Z))$ with the metric induced by $\exp^\perp$ from $X$ and the proof follows.\vspace {1mm}

{\sc Mean Curvature Inequalities.} All   of the above generalizes to Riemannian manifolds with {\it mean convex boundaries}, where the simplest instance of this is as follows.  

    \textbf {1.D.  Narrow Tunnel   Theorem.}  {\sl Let  $X\subset \mathbb R^n$, $n\geq 3$, be a smooth (not necessarily connected and   possibly infinite)    domain  with   
    $$mean.curv(\partial X)\geq n-2=mean.curv \left(S^{n-1}\left(R_+\right)\right)\mbox {   for $R_+=1+\frac {1}{n-2}$} .$$
If $n\leq 8$,    then the diameters of 
   all   {\it connected   components} of the subset 
   $$X_{-(1+\varepsilon)}=X\setminus U_{1+\varepsilon}(\partial X)=\{x\in X\}_{dist (x, \partial X) >1+\varepsilon } \subset X,
   \mbox { } \varepsilon>0,$$  are \it bounded by
$const_{n,\varepsilon}<\infty$}. (This  confirms  conjecture  10 in [Gr2019] for $n_1=1$.)

\vspace {1mm}

{\it Remark/Conjecture.}  Probably, the  connected  components of the subset $ X_{-1} \subset X$ are {\it bounded}.\footnote {Embarrassingly, I see no direct elementary proof of this even for $X_{-(1+\varepsilon)}$, where  $\varepsilon>0$.}

\vspace {1mm}
{\sc Plan of the Paper and Organisation  of the Proofs.} In  the  next (sub)section, we  discuss of our results  in the general context  
  of  the metric  scalar curvature geometry. \vspace{1mm}

In section 2  we introduce a (not quite)  numerical invariant on homology  of a  Riemannian manifold $X$, denoted  
$Sc_{\widetilde {sp}}^{\exists\exists\times}(h)$, $h\in H_m(X)$.

  We represent  $h$  in some cases by  {\it stable  $ \mu$-bubbles  $Y^m_{bbl}\subset X$}  with a    properly      tailored  function $\mu(x) $\footnote {Compare  with  sections
 5.3, 5.5 in [Gr2021].}  and thus    bound $Sc_{\widetilde {sp}}^{\exists\exists}(h)$ in terms of $Sc(X)$ in  manifolds $X$, where the universal covering  $ \tilde X$ is spin. \vspace{1mm}

{\sf (These $Y^m_{bbl}$ are constructed as in section 5 of [Gr2021] by
using the inductive dimension  descent scheme originated in  [SY1979$'$] and combined with  {\it symmetrization-by- $\mathbb T^{n-m}$-warping} from  [GL1983], also see  [SY2017],  following an  
idea   from [F-CS1980].
A presence of  possible stable  singularities in $ Y^m_{bbl} $   entails here  the probbaly  unnecessary   constraint $dim(X)\geq 8$.)}
\vspace{1mm}

In section 3 we introduce another invariant, denoted $Rad(h/S^N)$, which concerns area decreasing maps $X\to S^N$.
and show  that \vspace{1mm}  
   
   \hspace {25mm} $Rad(h/S^N)\leq \sqrt{m(m-1)/Sc_{\widetilde {sp}}^{\exists\exists\rtimes}(h)}. $\vspace{1mm}

{\sf (This is done essentially, but not quite,  as in   [GL1983] by relying   on the {\it relative index theorem for families of twisted Dirac operators}  over spin manifolds\footnote {This theorem generalizes to manifolds, where the universal coverings are spin,
 see \S\S$9\frac{1}{9}$, $ 9\frac {1}{8}$  in [Gr1996] and  section 10 in [Gr2019].} combined with the {\it "twisted" Lichnerowitz formula}, which is     sharpened  in the present case   by 
  {\it Llarull's inequality,}\footnote{Mario Llarull was Blaine Lawson's student who,  by Blaine's suggestion, algebraically evaluated   the bottom of the spectrum  of the 
   the  curvature  operator on the spinors on $S^n$, (formula 4.6 in [Ll1998]) and thus   obtained  the {\it optimal}   bound for the area dilations of maps from manifolds with $Sc\geq \sigma$ to spheres.}    
    compare with  section 3.4.3 in [Gr2021]).}

 \vspace{1mm}

 In section 4,  we show that the inequalities $sect.curv(X)\leq \kappa$
 and $ inj.rad(X)\geq \pi/\sqrt\kappa$ imply that $Rad(h/S^N)\geq 1/\sqrt\kappa$
 for $N\geq 2dim(X)$.
 
We combine this with the above and obtain the inequality \vspace{1mm}
  
 \hspace {33mm} $\kappa \geq Sc_{\widetilde {sp}}^{\exists\exists\rtimes}(h)/m(m-1) $ \vspace{1mm}
  
 \hspace {-6mm}   for  manifolds $X$  with the universal coverings  spin, which is a generalized quantitative version 
   of 
    the  following {\it non-existence theorem}  (see 13.8  in  [GL1983]).

 \vspace{1mm}

\OrnamentDiamondSolid \hspace {1mm} {\it Non-torsion} {\sf 
  homology classes\footnote{I am {\color {red!55!black} not certain} on what  is known in this respect about $l$-torsion, say for $l\neq 2$.} in complete manifolds  with non-positive

 sectional curvatures 
 are 
 not representable by continuous maps from oriented

  spin manifolds $X$ with positive scalar curvatures.}\footnote {The origin  of 
  the   argument used in [GL1983]  can be traced to the proof of the {\it Novikov conjecture} for spaces with $sect.surv\leq 0$ 
in [Mi1974].}
\vspace{1mm}
 
Finally, for $n\leq 8$,   we combine the inequality   $\kappa \geq Sc^{\exists\exists\rtimes}(h)/m(m-1) $ with a lower  bound on $Sc^{\exists\exists\rtimes}(h)$
obtained with $\mu$-bubbles  $Y^m_{bbl}$ in section 2,  \footnote {The    extra $\mathbb T^{n-m}$-warping factor associated with  $Y^m_{bbl}$
  is handled by    twisting the relevant (already twisted) Dirac operator with a {\it family of flat bundles} over  $\mathbb T^N$ or with a single {\it almost flat bundle} as in  [GL1980] and in  [Gr1996]. (See   [Gr2021] for  a pedestrian presentation of  these topics).} which is essential  to      gain the correct   $m(m-1)/n(n-1)$ factor  in our inequalities. 

 Thus we arrive at our \textbf {main theorem} 4.C  for distance  decreasing maps    from manifolds with lower bounds on their  
 {\it $ \mathbb T^\rtimes$-stabilized scalar curvatures} to manifolds with "large" metric  balls.

  \vspace {1mm}

In section 5 we say more on focal radii and discuss    lower bounds on curvatures of immersions. 

In section 6 we generalize  4.C   to manifolds with mean convex boundaries, for which we adapt the results from  [GS2000]  and    [Lo2021], which we  slightly improve along the lines of  [Ll1998]  and [Li2010].

In section 7 we   prove a   version of  1.B  for Riemannian foliations (defined in 1.1.B  below)  on 4-manifolds.

\vspace{1mm}

  Everywhere, we     indicate limitations of our results and arguments, (most pronounced for Riemannian  foliations with $Sc>0$) 
  and  try  to formulate   conceivable  generalizations and improvements.

\vspace {1mm}

\vspace {1mm}
\subsection {Corollaries,  Questions, Generalizations} 
\textbf {(a)}  Most known  inequalities  obtained 
 by  combining  (Dirac)  index theoretic   and ($\mu$-bubble) geometric measure theoretic 
 arguments  have been  generalized and   proved with the index theorems for {\it Callias type operators}  with no use of the geometric  measure theory.\footnote{See   [Ze2019]  [Ze2020] [Ce2020]  [CZ2021] [WXY2021].}

 However,  for all I know, there are no  such proofs  of 1.B, 1.C, 1.D   at the present moment  and the  constrain $n\leq 8$ remains intact.

 \vspace {1mm}

\textbf{(b)}  The  {\it $\square^{\exists\exists}(n,m=2,N)$-inequality} in 2.B of the next section and  
the ({\sf weakened}) {\it $\mathbb T^\rtimes$-stable 2d Bonnet-Myers diameter inequality for surfaces $Y$ with $Sc^\rtimes \geq 2$}   imply the following. 

$[\circ \bullet]_{m=2} $ {\it Pairs of Riemannian  metrics  $g$ and  $ \underline g\leq   g$ on  compact   $n$-manifolds  $X$  {\it  homeomorphic to} $S^2\times Z$, where  $ Z$ supports Riemannian metrics   $sect.curv\leq 0$, satisfy, for $n\leq 8$, the following:   
 $$ infSc(X,g)\cdot inj.rad (\tilde X,  \tilde {\underline g} )^2 \leq 2\pi^2 \frac { 2(n-1)}{ n},$$} 
where $\tilde X$ is the universal covering of $X$.\footnote{In fact,  
$ infSc(X,g)\cdot inj.rad (\tilde X,  \tilde {\underline g} )^2 <  2\pi^2 \frac { 2(n-1)}{ n} $, since, 
 by the proof of  [BMd], the inequality $Sc^\rtimes(Y) \geq 2$ implies that   $Y$  contains a pair of points with $dist(y_1,y_2) < 2\pi^2 \frac { 2(n-1)}{ n} $.
 
 {\color {red!55!black}Probbaly, } it is not hard to  show   that 
  $ infSc(X,g)\cdot inj.rad (\tilde X,  \tilde {\underline g} )^2 \leq   2\pi^2 \frac { 2(n-1)}{ n}-0.01 $, but  
  the corresponding 
 sharp inequality $ infSc(X,g)\cdot inj.rad (\tilde X,  \tilde {\underline g} )^2 \leq   { 2\pi^2} $
 remain {\color {red!55!black} problematic}.}

\vspace {1mm}
Furthermore, if $\inf Sc(g)/{\sf sup}sect.curv (\underline g)\geq n(n-1)$, then  {\it Zhu's  area inequality}\footnote {See
 [Zho2019],    [Zhu2020] and 2.8 in [Gr 2021]).}  implies    
1.B(i)  for $m=2$,  that is :
 $$ infSc(X,g)\cdot inj.rad (\tilde X,  \tilde {\underline g} )^2 \leq  2\pi^2.$$

Here, the Dirac operator is not used and the  spin issue doesn't arise.

  Besides,  granted 
a mild generalization  of    theorem 4.6   from [SY2017], the condition  $n\leq 8$ can be removed,\footnote{I can't claim this as a theorem, since I haven't  mastered 
the argument in [SY2017]. } while 1.B in the spin case may follow for all  $m$ and $n$   by adopting the techniques/ideas  of J. Lohkamp  [Loh2018]. (One has   no idea what to do with spin for $m\geq 4$.)

\textbf{(d)}   {\it \textbf {Contractibility $\lambda$-Radius.}} Given a function  $\lambda(r)\geq r$, $r\geq 0$,  define  $contr_\lambda rad(X)$ for a metric space $X$ as the {\it supremum of  the radii of  the balls $B _x(r)\subset X$,  which are contractible in
 the concentric ball $B_x(\lambda(r))\supset B_x(r)$  for all $x\in X$}.

Observe that  {\sf $contr_1 rad(X)\geq inj.rad(X) $} and that  
{\sf if $X$ is a contractible 
manifold with a cocompact isometry group, e.g. $X$  is  the universal   covering of  a compact aspherical manifold, then 
$contr_\lambda  rad(X)=\infty  $    for some function $\lambda(r) $.}

 Similarly to $[\circ \bullet]_{m=2} $ one shows that    if $n\leq 8$, then 
  
{\sl the metrics   {\sf  $g\geq  \underline g$ from 
 $[\circ \bullet]_{m=2} $ satisfy
 $$  contr_\lambda rad (\tilde X,  \tilde {\underline g} )  \leq \lambda \left ( 2\pi\sqrt {\frac { n-1}{ n\cdot infSc(X,g)\cdot}}\right),\leqno {[\circ_\lambda  \bullet]_{m=2}} $$
 for all function $\lambda(r).$}}
 (This generalizes $[\circ \bullet]_{m=2} $ for $\lambda(r)=r$.)

 Next let  $m=3$ and recall (see  3.10.1 in [Gr2021])  that the filling radii of 3-manifolds are  bounded in terms of  their $\mathbb T^\rtimes$-stabilized scalar curvature $Sc^\rtimes(X)\geq infSc(X)$ (defined in he  next section) as follows    
  $$fillrad(X)\leq \frac {6\pi}{\sqrt {Sc^\rtimes (X)}}.$$
 If $n=dim(X)\leq 8$, then the  $\square^{\exists\exists}(n,m=3,N)$-inequality  implies 
  (by an easy argument, see   (see [Gr1983]) that 
  $$ contr_\lambda rad (\tilde X,  \tilde {\underline g} )  \leq \Lambda \left ( \frac {6\pi}{\sqrt {Sc^\rtimes (X,g)}}\right),\leqno {[\circ_\Lambda  \bullet]_{m=3}} $$
 for all $\lambda=\lambda(r)$ and 
 $$\Lambda(r)=2\lambda(2\lambda(2\lambda(r)))$$
 (Desingularization  of $\mu$-bubles   needed for $n=dim(X)\geq 8$ seems here more difficult than that  for $m=2$.)

  \textbf {(e)} {\it \textbf {Stable Quasi-conjugate $\lambda$-Radii.}} There are  several
  candidates for a (bi)Lipschitz stable version of conjugate radius, where our definition below is, probably,  provisional.

  This  "radius"  of a Riemannian manifold $X=(X,g)$,  
   $$conj_\lambda rad(X)\geq conj.rad(X),  $$ 
   is the supremum of the numbers $r$, such that {\it there  
   exist  the following}:

  ($\circ$)      a smooth Riemannian manifold $\Theta$  of dimension $2n=2dim(X)$
   (instead of the space of  tangent $r$-balls, $TB(r)(X)=  \{ B_{0(x) }(r)\subset T_{x}(X)\}_{x\in X}$),  
   
  ($\circ$)     smooth maps
    $$\Pi: \Theta\to X, \mbox{ } \Upsilon:  \Theta\to X,  \mbox{ }  O: X\to \Theta,$$

   ($\circ$)     a {\it  homotopy} $\mathcal I_t :\Theta \to\Theta $, $ 0\leq t<\infty$, with the following properties.

  $\bullet_\Pi$   The map $\Pi$ is a smooth {\it submersion}, the fibers of which are denoted $\Theta_x=\Pi^{-1}(x)\subset \Theta$.

    $\bullet_\Upsilon$  The maps 
    $$\Upsilon_x=\Upsilon_{|\Theta_x} :\Theta_x\to X$$ are locally diffeomorphic on all  $\Theta_x$, $x\in X$.

 $\bullet_O$  The map $O:X\to\Theta$ (which corresponds to the zero  section 
 $\mathbf  0: X\to T(X)$) satisfies:  $O(x)\in \Theta_x=\Pi^{-1}(x) \subset \Theta $ for all $x\in X.$

 $\bullet_{\cal I}$  The homotopy starts from the identity map, $\mathcal I_0=id$, it preserves all  $\Theta_x$, i.e.      
 $\Theta_x\overset {\mathcal I_t} \to \Theta_x$ for all $t$ and  all  $x\in X$    and it fixes $O(X)$,  i.e. $O\circ \mathcal I_t(x)=x$ for all $x\in X$. 
    
  Moreover,  the maps $ \Theta\overset {\mathcal I_t} \to \Theta$ are  homeomorphisms\footnote{This bizarre  condition is need to justify 1.2.A below.}   for sufficiently large $t$,  and they   converges to $O$ for $t\to \infty$,  that is 
  $$\mathcal I_t(\theta)\to O(x)\mbox 
    { for  $  \theta\in \Theta_x $ and $t\to \infty$}.$$

  Now let us turn to   metric properties of  $O$ and $\mathcal I_t$ with
   respect to the induced metrics $g_x  =\Upsilon_x^\ast(g)$ on 
   $\Theta_x$.

  $\bullet_r$ The distances from  $O(x)\in \Theta_x$ to the boundary $\partial \Theta_x$ satisfy  
  $$dist (x,\partial\Theta_x)> r. 
  \footnote {Here we assume that $X$ is complete  with  no boundary  and this inequality means that the  $r$-ball in 
  $\Theta_x)$  around $O(x)$ is compact. }$$     
   
   $\bullet_\lambda$ For all $t\geq 0$, all $\rho>0$ and all $x\in X$, the map 
   
  $\mathcal I_t:\Theta_x\to\Theta_x$ {\it sends $\rho$-balls around   
  $O(x)\in  \Theta_x$ to the concentric $\lambda(r)$-balls.}

  
  \vspace{1mm}
  
  {1.1.A.} {\it \textbf {Observation.}}  The arguments in [CL2020]  and [Gr2020], as expounded     in  section  3.10.3 of  [Gr 2021], show that 
 
 {\it  for all functions $\lambda=\lambda(r)\geq r$,  there exists a positive  number $R=R_\lambda<\infty$,
such that all  compact Riemannian manifolds $X$ of dimensions $n\leq5$ with $Sc^\rtimes(X)\geq 0$  satisfy
 the following $\lambda$-version of Green's $[Sc|conjr]_n$:
 $$conj_\lambda rad(X)\leq {R_\lambda\over \sqrt {Sc^\rtimes (X)}}.
 \leqno {conj_{
 \leq 5}} $$}
 
  \vspace{1mm}
  
Notice that this does  not yield the corresponding inequality for the 
contractibility radius of (the universal covering of) $X$, since our definition  
 {\it doesn't  make} $conj_\lambda rad(X)\geq contr_\lambda rad(\tilde X);$ 
 However, the inequality $$contr_\lambda rad(\tilde X)\leq {R_\lambda\over \sqrt {Sc^\rtimes (X)}}.$$
 holds for $dim(X)\leq 5$   with the same $R_\lambda$  and  for the same reason as $conj_{
 \leq 5}$ does.

 {\sc Questions.}   All of the above,  including the $[Sc|conjr]_{n-k}$-conjecture, tells us precious little about    bounds on  "radii" $rad^\ast(X)$ of a Riemannian manifold  $X=(X,g) $ with $Sc(g)>1$, where such a bound must depend on   
  the topology and/or metric geometry of $X$ and 
 where   $rad^\ast $ may stand for  $inj.rad$, $conj.rad$, $inj.rad(\tilde X),$,   the  contractibility radius $contr_\lambda rad$ of $X$ or $\tilde X$  for some $\lambda(r)$
 or for $conj_\lambda rad$.

 Here are a few  specific  questions. \vspace {1mm} 
  
 {\it  Let $(X,g_0)$ be a compact Riemannian $n$-manifold with $Sc(g_0)>0$}.
  \vspace {1mm} 
  
  {\sf When can $g_0$  be $[Sc>0]$-{\it homotoped}\footnote{This means a homotopy  $g_t$  of 
  Riemannian metrics with $Sc(g_t)>0$, for all $t\in [0,1]$.} 
  to a metric $g_1$ with   
  $Sc(g_1)\geq n(n-1)=Sc(S^n)$  and $inj.rad(X,g_1)  \geq \rho$  for a given $\rho>0$?}
  
(I)   {\it For instance}, -- {\sf this seems unlikely but not impossible } --  does such  $[Sc>0]$-homotopy $g_t$  exist for {\it all} $\rho<  \pi$ and all {\it simply connected }manifolds  $X$?
 
 (II) Or, does the inequality $Sc (g)\geq n(n-1$  for a metric  $g$ on  $X$ imply that 
  $conj.rad (X,g_1)\leq \pi-\varepsilon_n$ for some universal
  $\varepsilon_n>0$, e.g.  for  $\varepsilon_n=1/10^n$, unless
  $g$ is $[Sc>0]$-homotopic to (is closed to?) a metric with constant sectional curvature on $X$?

(III) What are {\it $Scalrad^\ast$-extremal} manifolds $(X,g_0)$, i.e.  such
 that  
 {\it $$\inf Sc(g_0)\cdot rad^\ast(X,g_0)^2\geq \inf Sc(g)\cdot rad^\ast(X,g)^2$$
 for all Riemannian metrics $g$ on $X$? }
 
 For instance, are {compact irreducible  symmetric spaces} 
   $Scalrad^\ast$-extremal,  say for $rad^\ast=inj.rad$?\footnote{This may be also true for some reducible spaces such as $S^2\times S^2$, for instance.} 
   
 (IV) How much does  a bound on the diameter of $X$ contribute to the inequality $ rad^\ast(X)\geq \rho  $  for manifolds with $Sc(X) \geq 1 $?
 
 For instance,  what can be said about the {\it  topological  complexity of an $X$ which admits a metric $g$, such that 
 
 $$Sc(g)\geq 1,  \mbox { }   inj.rad(X,g)\geq 0.01,   \mbox  { } diam(X,g) \leq 1?$$}
 
 Do these inequalities limit  some  complexity of the $[Sc>0]$-homotopy class of $g$?
  
  {\it Remark.}  If we drop  the inequality $Sc(g_1)\geq 1$,  the remaining two, which amounts to ${diam\over inj.rad}\leq 100$, don't much restrict the topology of $X$.
  
For instance all     surfaces $X$ admit Riemannian metrics with  ${diam\over inj.rad}\leq 4$  obtained with ramified coverings  $X\to \mathbb T^2$ 
and/or $X\to S^2$  (see [Ba1996]). 

Similarly, by ramifying $X^n\to S^n$ one can produce topologically complicated manifolds with  small  ratios ${diam\over inj.rad}$.

In fact, using  {\it Alexander's   braid-knot theorem +  Hilden's 3-fold branched covering  theorem} (applied to $S^3$ preliminary very strongly ramified over two linked Hopf circles)  one can show that 

 {\sl  {\it all}  compact 3-manifolds  admit Riemannian metrics  with  ${diam\over inj.rad}\leq 100$.}

I am not certain if this is true for all $n$-manifolds.

\textbf {1.1.B.  Riemannian Foliations.}  Let   
$$\mathcal X= (\mathcal X,g) = (Q,\mathcal X,  g)$$  
be a smooth    Riemannian $n$-foliation, that is a smooth manifold $Q$ foliated  into smooth $n$-dimensional manifolds $X$, and $g $ is a smooth positive definite   quadratic form on the (sub)bundle
 $T(\mathcal X)\subset T(Q)$ tangent to the leaves.

Let $X_q\subset Q$, $  q\in Q$  denotes  the leaf passing trough the point  $q$ 
 and let $Sc(\mathcal X, q)$, $conj.rad_q(\mathcal X)$, etc.  be the corresponding invariants   of $X_q$  at $q$ defined in section 1, where the corresponding $inf$-invariant refer to $\inf_{q\in Q}.$  

Observe that Green's integration argument applies to Riemannian foliations with {\it transversal measures} on compact manifolds\footnote{Here one needs only  transversal  continuity of $\mathcal X$ and $g$.}  and that 
{ \color  {red!65!black} $[Sc|conjr]_{n-k}$-Conjecture 1.A} makes sense for foliations.

\vspace {1mm}


\section { $T^\rtimes$-Stabilized Curvatures $Sc^\rtimes$,  $Sc^{\exists\rtimes} $,  $Sc^{\exists\exists\rtimes}$ and  
Multispreads $\tilde \square^\perp(h)$ on Homology} 
 
{\it $\mathbb T^\rtimes$-Extensions}.   A  {\sf "warped"  $\mathbb  T^N$-extension}, $N=0,1,....$, where $\mathbb T^N$ is the $N$-torus, of  a Riemannian manifold $Y=(Y,g)$, where  $Y$ may have a boundary,    is a Riemannian manifold 
 denoted $Y_N^\rtimes=Y\rtimes \mathbb T^N$, that is the product $Y\times \mathbb T^N$
with  a Riemannian metric $g_N^\rtimes $ of the following form:
$$g_N^\rtimes=dy^2+\sum_{i=1}^N\varphi_i^2dt_i^2$$  
for some  smooth positive functions  
 $\varphi_i(y)\geq 0$  on $Y$, which are strictly positive (>0) in the interior $Y_N\setminus \partial Y$ of $Y$.

Clearly,  this $g_N^\rtimes$ is invariant under the obvious  action of $\mathbb T^N$  on $Y_N^\rtimes  $ and  
 $\varphi_i(y)$ are  equal to  the  $g^\rtimes$-lengths of the $\mathbb T_i^1$-orbits of the points  
 $y^\rtimes = (y,t) \in Y^\rtimes$ under  this action.

{\it $\mathbb T^\rtimes$-Stabilized Scalar Curvature(s)}.  Let  
us agree that  the inequality 
$$Sc_N^{\exists\rtimes} (Y,y)>\sigma(y) $$
 for a  given  function  $\sigma(y)$ on $Y$ signifies that 

{\it  there exists  a $ \mathbb T^N$-extension $Y_N^\rtimes$ of $Y$
the scalar curvature of which satisfies 
$$Sc(Y_N^\rtimes,  y)> \sigma(y),$$}
where this scalar curvature is a function $Y=Y_N^\rtimes/\mathbb T^N$, since the
 metric $g_N^\rtimes $ is  $ \mathbb T^N$-invariant.\footnote { In fact, 
 $Sc\left(g+\sum_{i=1}^N\varphi_i^2dt_i^2\right)= Sc(g) -2\Delta\Phi-||\nabla \Phi||^2 -\sum_i{||\nabla\log \varphi_i||^2}
 $
 for $\Phi(y)=\log(\varphi_1(y)\cdot ...\cdot \varphi(y)_N)$, see
 [GL1983], [SY2017], [Zhu2019], or  compute yourself. }

Clearly 
$$ Sc(Y)\leq Sc_1^{\exists\rtimes}(Y)\leq...\leq  Sc_N^{\exists\rtimes}(Y)\leq ...\mbox {  .}$$ 

We agree that $Sc^{\exists\rtimes}(Y)$ stands for $Sc_\infty^{\exists\rtimes}(Y)=\sup_N Sc_N^{\exists\rtimes}(Y)$, 
i.e. $Sc_N^{\exists\rtimes}$  with 
 an arbitrarily  large $N$ and observe that  $Sc^{\exists\rtimes}(Y)$   (and even $Sc_1^{\exists\rtimes}(Y)$)
 of manifolds $Y$ with boundaries (see below)  can be  significantly greater than $Sc(X).$
Yet, most known geometric properties of manifolds with $Sc> \sigma$ are also enjoyed  by those with 
 $Sc^{\exists\rtimes}> \sigma$.\footnote{The proofs of the  bounds on the 2-waists  of 3-manifolds 
 $Y$ with $Sc(Y)\geq \sigma>0$ (see [MN2011],  section 3.10 in [Gr2021] and LM2021]) do not apply, 
 at least not immediately, to $Y$ with $Sc^\rtimes(Y)\geq \sigma$. 
 
 Also the $\mathbb T^\rtimes$-stabilization  of the known   geometric 4d-inequalities  obtained  with the Seiberg-Witten equations (see [LB2021])  remains {\color{red!55!black} problematic}.}
 
\vspace {1mm}

To have a single non-ambiguous number, let 
{\it $$Sc^{\rtimes} (Y)$$   be  the supremum} of the numbers $\sigma${, such that  
 $Sc^{\exists\rtimes} (Y)>\sigma$.\vspace {1mm}

{\it \textbf { Examples.}} It is not  hard to show the following (see [Gr2023]).  
 
  (a) {\it Rectangular solids} satisfy:
$$Sc^\rtimes\left( \bigtimes_{i=1}^n[-a_i,b_i]\right)=4\lambda_1\left( \bigtimes_{i=1}^n[-a_i,b_i]\right)=\sum_{i=1}^n \frac {4\pi^2}{(b_i-a_i)^2}.$$
where this $\lambda_1$ is the first eigenvalue of the Laplace operator with the Dirichlet boundary condition in $\bigtimes_{i=1}^n[-a_i,b_i].$

(b) Unit hemispheres satisfy 
$$Sc^\rtimes(S^n_+)=n(n+3)=Sc (S^n_+)+4\lambda_1(S^n_+).$$ 

(c)  Unit balls  satisfy
$$Sc^\rtimes (B^n)=4\lambda_1(B^n)> n^2+5n+n^{4/3}.$$

{\it \textbf { Definition  of} $\mathbf { Sc^{\exists\exists}.}$} Given a homology class in a Riemannian manifold $X$,
$$h\in H_m(X),$$ 
write
$$Sc_N^{\exists\exists\rtimes}(h)=Sc_N^{\exists\exists\rtimes}(h)_{dist}>\psi \mbox {  for a function $\psi=\psi(x)$ on $X$},$$
if {\it there exits}  a compact oriented (a priori disconnected)  Riemannian $k$-manifold $Y$ and  a 1-Lipschitz (i.e. distance 
non-increasing)   map 
$$f:Y\to X, \mbox { such that
 $f_\ast[Y]=h$},$$ 
  where $[Y]\in H_k(Y)$ is the fundamental class of $Y$, 
  and such that  the $\mathbb T^N $-stabilized scalar curvature of $Y$ is bounded from below by $\psi$, i.e. 
$$Sc_N^{\exists\rtimes}(Y)> \psi \circ f$$
for the composed function $\psi \circ f(y) = \psi(f(y))$.\vspace {1mm}
 
 $\mathbf  {[Sc^{\exists\exists\rtimes }(X)]}$. Here $Sc^{\exists\exists\rtimes }(X)=Sc^{\exists\exists\rtimes }[X]$ for the fundamental homology class 
 $[X]\in H_{dim(X)}(X)$ of an oriented manifold $X$.  

Observe that
$$Sc^{\exists\exists\rtimes }(X)\geq Sc^{\rtimes }(X),$$
where this inequality is strict, for instance, for products $X=X_1\times X_2$, where $X_1$ has constant scalar curvature 
$\sigma_1>0$ 
and $Sc(X_2)<0$.

On the other  hand the  equality  {\color {red!55!black} is expected} for many manifolds with positive (sectional, Ricci?) 
curvatures, but at the present moment the proof is available only for spheres $S^3$, $S^4$ and closely related
 manifolds of dimensions $3$ and $4$. Yet higher dimensional  equalities  of this kind 
 e.g.  for compact symmetric spaces with non-zero Euler characteristics, are available (see{GS2000] and section 8.1) for the {\it spinor versions 
 of $Sc^{\exists\exists\rtimes}$}
 defined below.

 {\it Open manifolds and Homology with Infinite Support.}
 The above definition of $Sc^{\exists\exists\rtimes }$ naturally generalizes  to homology classes  with {\it infinite supports} in non-compact  manifolds, denoted 
 $$h\in H_m(X,\partial_\infty X),$$  
 where $\partial_\infty X$ stands for the  complement of an unspecified  "arbitrarily large" compact subset in $X$.
 
 In other words,  $H_m(X,\partial_\infty X)$ is   the projective  limit of $H_m(X, X_i)$  over all  $X_i\subset X$ with {\it compact complements.}

Here the manifolds $Y$ may be non-compact and the maps $f:Y\to X$ are required to be proper,  that is $f$ sends $\partial_\infty Y\to \partial_\infty X$, which means that the $f$-pullbacks of compact subsets in $X$ are compact.\vspace {1mm}
 
  
  

\vspace {1mm}

{$\mathbb \exists\exists^\rtimes$}-\textbf {Spin.}  One defines two spin versions of   $Sc^{\exists\exists\rtimes }$:  
 
 $Sc^{\exists\exists\rtimes }_{{sp}}(h)$, where the manifolds $Y$ are  {\it spin},
  
 $Sc^{\exists\exists\rtimes }_{\widetilde {sp}}(h)$,  where  the {\it universal coverings} $\tilde Y$ of $Y$ are spin.

  Clearly,
  $$Sc_{sp}^{\exists\exists\rtimes }(h)\leq Sc^{\exists\exists\rtimes }_{\widetilde {sp}}(h)\leq Sc^{\exists\exists\rtimes }(h).\footnote{ {\color {red!55!black}Conceivably,} some  nonzero multiples of all   homology classes $h$  satisfy $Sc_{sp}^{\exists\exists\rtimes }(i\cdot h)= Sc^{\exists\exists\rtimes }(i\cdot h)$.  }$$

\vspace {1mm}

{\it Homological Pullbacks.} Let  $f:V\to W$, where $dim (W)-dim (V)=k$,  be  a proper  map between orientable manifolds,  and recall  that  the cohomology homomorphism 
$f^\ast: H^l(W) \to H^l (V)$ composed  with the Poincare duality isomorphisms defines the "pullback homomorphism" $H_m(W)\to H_{m+k}(V)$ for $m=dim(V)-l$,   which, for generic smooth maps, is implemented by taking $f$-pullbacks of $m$-cycles in $W$.

Also,  such a homomorphism is defined for topological fibrations and   topological submersions  $V\to W$ with manifold fibers. 
\footnote {This is explained in purely geometric terms   in [Gr 2014].}

{\it \textbf {Multi-Spreads.}} (Compare with  section 7.1 in [Gr2021].)
Let $V$ be an $n$-dimensional Riemannian manifold (possibly non-compact)  with a boundary,   and let $ h\in  H_m( V, \partial_\infty  V)$ be a  homology class (with  infinite, i.e. non-compact, support if  $\partial_\infty$ is  non-empty,  i.e.  $V$ is non-compact).

 {\it The $\square^\perp$-spread} 
$\square ( h)$  
is  the {\it supremum   of the  numbers} $d\geq 0$, for which  {\it there exists}

{\sf a continuous  boundary-to-boundary map  from $V$ to the $k$-cube for $k=n-m$,
$\psi=(\psi_1,...\psi_i,...\psi_k): V\to  \square =[-1,1]^k$,   
$\psi_i:  V\to [-1,1]$,
such that the following two conditions are satisfied.

(1)   The  class $h$ is equal  to the  $\psi$-pullback of the generator of $h_0\in H_0([-1,1]^k)$,  that is, 
 if $\psi$ is smooth, then  $h$  is represented by
 the $\psi$-pullback of a generic point $p$ in the interior of the cube.
  
(2)  {\sf The distances  between the pullbacks of the opposite faces in the cube $[-1,1]^k$, 
$$d_i=dist_V(\psi_i^{-1}(-1),\psi_i^{-1}(1)), i=1,...,k,$$
are bounded from below by the following inequality 
$$ \left (\frac {1}{k}\sum_{i=1}^k \frac {1}{d_i^2}\right)^{-\frac {1}{2}}\geq d,$$
that is 
$$\frac {1}{k}\sum_{i=1}^k \frac {1}{d_i^2}\leq \frac {1}{d^2},$$}
(e.g. if $d_i\geq d$ for all $i=1,...,k$).}

(Equivalently, one could require the functions  $\psi_i: V\to [-1,1]$ to be 
$\frac{2}{d_i}$-Lipschitz.)

\vspace{1mm}

Next define $\tilde\square^\perp(h)\geq \square^\perp(h) $ of a homology class $h\in H_m(X,\partial_\infty X)$
for an $n$-dimensional Riemannian manifold $X$ (possibly non-compact)  with a boundary, 
  as the supremum of the  numbers $\tilde d\geq 0$, 
such there exist 
 
 (i)  a Riemannian manifold  $\tilde V$ of dimension $n=dim (X)$   
 
 (ii)  a homology class $\tilde h\in H_{m}(\tilde V,\partial_\infty \tilde V)$ with $\square^\perp( \tilde h)\geq d$, 
 
(iii)  a proper locally isometric map\footnote {This map {\it does'n have} to send $\partial \tilde V\to\partial X$.} $\phi:\tilde V\to X$,  such that  the induced homology homomorphism 
$\phi_\ast :H_{m}(\tilde V;\partial _\infty \tilde V) \to H_{m}( X; \partial_\infty X )$
sends $\tilde h$ to $h$, in writing:
$\phi_\ast(\tilde h)=h.$

 \vspace{1mm}

\textbf {2.A.} {\it Enlargeable Manifolds  and $sect.curv\leq 0$.}  A Riemannian manifold $Z$ is called enlargeable if  there 
exist zero dimensinal homology classes $h_i\in H_0(Z)$ with $\tilde\square^\perp(h_i)\geq i$,  $i=1,2,....,$. 

For instance, complete manifolds with non-positive sectional curvatures are  enlargeable, since their universal coverings  contain {\it arbitrarily  large} 
compact domains $\tilde V$, i. e. with arbitrarily large $\perp$-multispreads of their zero-dimensional homology.

It follows that if  $X\underset {homeo}= Y\times Z$, where $Y$ is an oriented $m$-manifold and  where $Z$ is 
a compact manifold, which admits a metric with non-positive 
sectional curvature, e.g. it is the $(n-m)$-torus, then  
the $\tilde \square^\perp$-spread 
(unlike the  $ \square^\perp$-spread) of
the {\it fundamental homology class}  $[Y]\in H_m(X,\partial_\infty X)$ of $Y=Y\times z_0
\subset  X$,   is {\it infinite.}

 \vspace{1mm}

  \textbf {2.B. $\square^{\exists\exists}(n,m,N)$-Inequality.}  {\sf Let $X=(X,g)$ be 
  an $n$-dimensional orientable  Riemannian manifold   with a boundary and let $h\in H_m(X,\partial_\infty X)$.}

If $n=dim(X)\leq 8$, then
      $$Sc_{N+k} ^{\exists\exists\rtimes} (h)\geq Sc_N^{\exists\rtimes}(X)
-  \frac {n+N-1}{n+N }\cdot\frac {4k\pi^2}{\tilde \square^\perp (h)^2}, \mbox { } k=n-m, \leqno{\color 
{blue}[Sc^{\exists\exists}]}$$} 
that  signifies  the implication 
$$Sc_{N}^{\exists\rtimes}(X)> \psi \implies Sc_{N+k}^{\exists\exists\rtimes} (h)> \psi$$ 
   for all continuous functions $\psi=\psi(x)$ on $X$.\vspace {1mm}
 
 Furthermore, 
 
 {\sl if $X$ is spin, then 
 $$Sc_{N+k,sp}^{\exists\exists\rtimes} (h)\geq Sc^{\exists\rtimes}_N(X)
-  \frac {n+N-1}{n+N }\cdot\frac {4k\pi^2}{\tilde \square^\perp (h)^2},\leqno{\color {blue}[Sc_{sp}^{\exists\exists}]}$$
and if the universal covering $\tilde X$ of $X$ is spin, then 
$$Sc_{N+k,\widetilde {sp}}^{\exists\exists\rtimes} (h)\geq  Sc^{\exists\rtimes}_N(X)
-  \frac {n+N-1}{n+N }\cdot\frac {4k\pi^2}{\tilde \square^\perp (h)^2}.\leqno{\color 
{blue}[Sc_{\widetilde {sp}}^{\exists\exists}]}$$}
 
 {\it Proof.}  If $k=1$ this follows from the  equivariant separation theorem applied to $X\times \mathbb T^ N$  with a 
 metric  $g^\rtimes$  and the general case follows by induction on $k$.\footnote {See 5.4 in [Gr2021]   also [Gr2021$'$]  and compare with [Ri2020] 
 and [GZ2021], where similar arguments are used}
  \vspace {1mm}
 
(See the following section for immediate corollaries of this and for open questions.)\vspace {1mm}

It is  commonly   {\color {red!55!black} believed} that a presence of singularities in minimal hypersurfaces and/or in stable $\mu$-bubbles  for $n\geq 8$ is {\it non-consequential  for scalar curvature}, which is confirmed  by  Natan Smale's generic regularity for $n=8$ 
[Sm2003]. 
  
In the present case, this belief turns into the  following.

 \textbf {2.C. $\square^{\exists\exists}(n>8)$-{\color {red!55!black}Conjecture.}} {\sf The   inequalities 
{\color {blue}$[Sc^{\exists\exists}]$}, {\color {blue}$[Sc_{sp} ^{\exists\exists}]$} and {\color {blue}$[Sc_{\widetilde{sp}}^{\exists\exists}]$}
 hold  true for all $m$, $n$, and $N$.}

(It is {\color {red!55!black} conceivable}  that   the inequality  {\color {blue}$[Sc_{sp}^{\exists\exists}]$}  holds for non-zero multiples $i\cdot h\in H_m(X)$ in  {\it non-spin} manifolds  $X$.)

 \textbf {2.D.   $\exists${\large $\rtimes $}-Stabilization of Riemannian Foliations.}  Given a Riemannian foliation  
  $(Q, \mathcal Y, g)$ 
with $m$-dimensional  leaves $Y$ (see 1.1.B), define its
 $\mathbb T^\rtimes$-extension}  
 $$(Q^\rtimes =Q\times \mathbb T^N, \mathcal Y^\rtimes, g^\rtimes), $$ 
by foliating  $Q^\rtimes$  into $Y^\rtimes= Y\times \mathbb T^N$  
 and taking $g^\rtimes=g+\sum_{i=1}^N\varphi_i^2dt_i^2$
  for   smooth positive functions  
 $\varphi_i(q)\geq 0$  on $Q$.  
 
  Accordingly, we introduce $Sc^{\exists\rtimes}(\mathcal Y)$, where the sentence
 
\hspace {27mm}  {\sf  "$Sc^{\exists\rtimes}(\mathcal Y)$ has property $P$}"
  
    means that {\it  there exist functions $\varphi_i$,  such that
    
  \hspace {20mm}  {\it  the corresponding   $\mathcal Y^\rtimes $ has property  $P$.} }
  
Conceivably, one may also define $Sc^{\exists \exists\rtimes}(h)$, $h\in  H_M(R)$,  for Riemannian manifolds $R$  by considering oriented $M$-dimensional manifolds $Q$ foliated by $m$-dimensional Riemannian $Y$  and maps $f:Q\to R$,  which are distance decreasing on all leaves  $Y$
with respect to the metrics $Sc^{\exists \rtimes}\cdot g$ on $Y$.

Then one may formulate a counterpart of $\square^{\exists \exists}$-inequality 2.B  for manifolds $R$ foliated into $n$-dimensional leaves  
for $dim(R)-n=M-m$. 

But there is little evidence for such a generalization of 2.B, where a more realistic conjecture should allow foliations with singular leaves $Y$.  (See  section 7.)


\section{Maps to Spheres and to Sphere Bundles }

Let $X=(X,g)$ and $S= (S,g_S)$ be   Riemannian manifolds, let $\mathcal P\subset C(X,S))$  
be a subset in  the space of continuous maps 
$X\to S$,  let $\Upsilon:X\times \mathcal P\to S$ be the evaluation map $\Upsilon:(x,a)\mapsto a(x)$.  
Let $B$ be an aspherical  topological space and $\beta:X\times \mathcal P\to B$ a continuous map, e.g. $\beta$ is the composition of   the coordinate projection $X\times \mathcal P\to X$ a map
$\underline \beta:X\to B$  and the classifiyng map  map
$\underline \beta:X\to B=B(\Gamma)$,  where $\Gamma $ is the fundamental  group  $\pi_1(X)$.

\vspace{1mm}

What we do  in this section is motivated by the following.

 \textbf {3.A. Aspherical Parametric Mapping  Conjecture.} {\sf Let $X$ be a   complete  $n$-manifold,  let $S$ be the unit $K$-sphere $S^K $, $K\geq 1$,  let  $[S^K]^\ast \in H^K(S^K)=\mathbb Z$ be the fundamental cohomology class and 
let  $\alpha \in H^{k}(B)$.}

{\sl If    all  $f\in\mathcal P $  
 are   smooth maps $f:X\to S^K$,  which are constant at infinity and which  are {\it strictly area decreasing} on all  smooth surfaces in $X$
with respect to the Riemannian metrics $Sc(g)\cdot g$ on $X$   and $ m(m-1)g_{S^K}$ on $S^K$ for $m=n-k$,  
then the cup product 
$$\Upsilon^\ast [S^K]^\ast \smile \beta^\ast(\alpha) \in H^{K+k}(X\times \mathcal P, \partial_\infty X\times P)$$ 
vanishes on 
$$H_n(X)\otimes H_{K-m}(\mathcal P)\subset H_{K+k}(X\times \mathcal P).$$}
(Recall that $m=n-k$  and that  $\partial_\infty X$  stands  for the complement  of a large unspecified  compact   subset in $X$.)

{\it Remarks} 
(i$_A$) Theorem 3.E  below  confirms  this conjecture, where $\alpha$ is the fundamental  cohomology class of an {\it enlargeable 
$k$-manifold},  provided that  $dim(X)=n\leq 8$ and the universal covering of $X$ is spin.

 (ii$_A$) This conjecture makes  sense for singular spaces $X$, where one can define  scalar curvature, e.g. for Alexandrov spaces  with lower curvature bounds.

(iii$_A$) An ultimate form of this conjecture must apply to the space/category  of germs of $n$-submanifolds (or singular spaces) in $S^K$ with induced metrics having $Sc\geq  \sigma>0$ (e.g. $\sigma=m(m-1)$).

 The following results (corresponding to  $m=n$)  gives an idea of what may be  expected in this  regard.      

  \textbf {3.B.  Su-Wang-Zhang Foliated Area Contraction Theorem.}  {\sf  Let $\mathcal X=( \mathcal X, G)$  be a complete oriented  Riemannian $K$-manifold 
 smoothly foliated into $m$-manifolds $\tilde X\subset \mathcal X$ and   let $f:\mathcal X\to S^K$ be a smooth map locally constant at infinity.}
 
{\it  If the leaf-wise scalar curvature of  $\mathcal X$, denoted $Sc^\sim(\mathcal X)$    is  bounded from below by 
$Sc^\sim(\mathcal X) >m(m-1)$ and if the map $f$ is area decreasing an all smooth surfaces contained in the leaves, then, provided either the manifold  $\mathcal X$ is spin or the tangent bundle to the leaves is spin, 
 the map $f$ has zero degree,
$$deg(f)=0.$$}

{\it Remarks.} (i$_B$)   The inequality $Sc^\sim(X)>m(m-1)$ is required  in [SWZ2021]  only on the support of the 
differential of $f$, while  everywhere else
   $Sc^\sim$  must be non-negative.

(ii$_B$) The conclusion $deg(f)=0$   holds  (unless I misunderstand something in [SWZ2021])  if  the map $f$ is leaf-wise area decreasing for the metrics 
$Sc^\sim \cdot G$ in $\mathcal X$ and $m(m-1)g_{S^K}$ in $S^K$. 

(iii$_B$) Probbaly, (I haven't properly  checked this) the spin condition can be relaxed to the spin of the lifts of the tangent bundles (of the manifold $\mathcal X$  and of the  leaves)  to the universal covering of $\mathcal X$.\vspace {1mm}

\textbf {Preparion for 3.C.} Let  $ P$ and $ Q$ be smooth manifolds 
possibly with  boundaries
and $\Psi :  Q\to P$ be  a smooth {\it submersion},\footnote  {The differential $d\Psi:T(Q)\to T(P)$  has everywhere $ rank=dim(P)$} the fibers of which,  denoted  
$$\tilde X_p=\Psi^{-1} (p) \subset Q, \mbox { }p\in P,$$
are thought of as a family of manifolds parametrized by $P$.

Let $\mathcal X$  be  the foliations of $Q$ into (the connected components of) the  fibers $\tilde X_p$ and let  $ g$  be a Riemannian metric on  $\mathcal X$. 

(In  the present case, this is  
 a  smooth family of  Riemannian metrics $g_p$  in  $\tilde X_p$.) 


  Let $S^\bullet \to P$ be a $K$-dimensional sphere bundle, let  $\delta: P\to S^\bullet$  be
 a  section of this bundle and let  $P_\bullet$ be the image of the opposite section,
 $$P_\bullet=-\delta(P)\subset S^\bullet.$$
  
  Let $\theta^\bullet\in H^K(S^\bullet, P_\bullet)  $ be
   the {\it Thom cohomology  class} of the $K$-ball   bundle $S^\bullet\setminus  P_\bullet  \to P$, that is the  Poincar\'e  dual to 
  $$\delta_\ast[P]\in H _d(S^\bullet,\partial S^\bullet),$$
    where $d=dim( P)$  and 
  $ [P] \in H_d (P, \partial P)$ is the fundamental class of $P$.
  
 Let   

 $$ H_{i/j}( Q, Q\setminus Q_0)= H_{i/j\downarrow\Psi}( Q, Q\setminus Q_0),\mbox { } Q_0\subset Q,$$
be the group of relative  homology classes of   possibly infinite relative   $i$-cycles 
$C\subset Q$, such that the  images  $\Psi(C)\subset P$ are {\it precompact}  and such that 
$dim(\Psi(C))\leq j$
and  $dim (\Psi(\partial C))\leq j-1$.

 Let $$\mathcal F :Q\to S^\bullet $$
  be a  smooth {\it fiber bundle} map,\footnote {$p$-Fibers of $Q$ go to  $p$-fibers of $S^\bullet$ for all $p\in P$,} 
 such that 
 
{\sf  (i)  the pullback  $ \Psi^{-1}(\partial  P) \subset Q$ is sent by $\mathcal F$ to $P_\bullet\subset S^\bullet$,
$$ \mathcal F( \Psi^{-1}(\partial  P)) \subset P_\bullet;$$

 (ii) the complements of "uniformly  controllably large" compact subsets $R_p$  in the fibers $X_p\subset Q $ are also sent to $P_\bullet$,  i.e. 
 
 {\sl  there is a closed subset, $R\subset Q$ such that the restriction
  $\Psi_{|R}:R\to P$ is a {proper} map and $\mathcal F(Q\setminus R)\subset P_\bullet$.}}


Let 

$$\mathcal F^\ast(\theta^\bullet)\in H^K(Q, \mathcal F^{-1}(P_\bullet))$$
be the  (relative)   cohomology class   induced  by $\mathcal F$  from $\theta^\bullet\in H^K(S^\bullet, P_\bullet)$.

{\it \textbf  {Spherical Metric $\gamma$.}} Let $\gamma=\{\gamma_p\}$,    be a smooth family of Riemannian metrics  in   the
  $K$-spheres 
  $S^\bullet_p\subset S^\bullet$  with  constant sectional  curvatures  $\kappa_p$.

\textbf {3.C. Prelimenary   Area Contraction   Inequality. } {\sf
 Let the  above  map $\mathcal F$ be 
 {\it fiberwise area decreasing}  with respect to the metrics 
 $Sc^{\exists\rtimes} (\mathcal X) \cdot g$ 
  in $Q$\footnote {Recall that $Sc(g^{\rtimes})$  actually is a function on $Q$, see 2.D.} and $ n(n-1)  \kappa\gamma=\{\kappa_p \cdot  n(n-1) \gamma_p\}$ in $S^\bullet$,
i.e.  the   $\mathcal F$-images of all smooth surfaces $E\subset \tilde X_p$, $p\in P$, satisfy
$$n(n-1)\kappa_p\cdot area_{\gamma_p}(\mathcal F(E))  <  \int_E Sc^{\exists\rtimes}  (\tilde X_p)(x_p)  d_{g_p}x_p.\leqno{[area<]_{\exists\rtimes}}$$
{\sf Let the  fibers $\tilde X_p$ be   manifolds {\it without boundaries} and the  metrics $g_p$ on $\tilde X_p$ be {\it complete}.

If   the sphere bundle $S^\bullet\to P$ and the   tangent bundle $T(\mathcal X)\to Q$ are {\it spin},
then  the cohomology class $\mathcal F^\ast(\theta^\bullet)$
 {\it vanishes  on 
 $H_{K/K-n}(Q, \Psi^{-1}(\partial P))$.}}}

\vspace {1mm} 

{{\it About the Proof.}  If the bundle   is trivial, then 3.C for area decreasing maps 
$(\mathcal X, Sc (\mathcal X) \cdot g)\to  
  Q) $ (instead of $ [area<]_{\exists\rtimes}$)   follows  from 3.B which need o be  augmented with (ii$_B$);  this, as explained
   in ??? below, yields 3.C for     $Sc^{\exists\rtimes}$ as well. 

In fact, since the leaves $\tilde X_p$  are closed, the proof of 3.C is much simpler than that of 3.B in [SWZ2021]) 
(where non-trivial holonomy creates the major problem) and it goes along the usual lines, roughly,  as follows.

 Let $\mathbb S^\mp(S_p^\bullet)\to P$ be the Clifford bundle 
of spinors in the 
fibers   $S^\bullet_p=S^K$, let $\mathcal T^\rtimes\to  \mathbb T^N$ be an
 {\it almost flat bundle with non zero top Chern class and all other classes zero} and 
 let $\mathcal D^\rtimes _p$, $p\in P$,  be the Dirac operators on the fibers $\tilde X^\rtimes _p=\tilde X_p\times \mathbb T^N$ with coefficients in  $\mathbb S^\mp(S_p^\bullet)\otimes \mathcal T^\rtimes \to P$,
   $$\mathcal D^\rtimes _p: C^\infty (\mathbb S^\mp(X^\rtimes _p)\otimes \mathbb S_{\mathcal F(p)}^\mp(S_p^\bullet)\otimes \mathcal T^\rtimes)\circlearrowleft.$$
   Then, by  Llarull's trace evaluation   (4.6 in [Ll1998]) of the curvature term in  the twisted Lichnerowitz formula 
   (see [GL1980], [GL1983] and   4.16 and  8.17 in [LM1989]), the inequality $[area<]_{\exists\rtimes}$ implies that  
   the operators  $\mathcal D^\rtimes _p$ are positive. (Zhang's  Dirac deformation argument [Zha2020]  allows 
    $\inf_xSc(X,x)=0$ at infinity, but this is  unneeded for  our applications.)

  Hence, the $K(P)$-theoretic  relative $\mp$index (see \S13 in [GL1983] and references therein)   is zero. 
   Then, if $K$, $N$ and  $dim(\tilde X_p)$  are even, the relative version of the Atiyah-Singer  index theorem for families  [GL1983]
  shows that 
   $$\mathcal F^\ast(\theta^\bullet)(h)=0\mbox {  for all }
 h\in H_{K/K-n}(Q, \Psi^{-1}(\partial P)),$$
   
   Finely,  a  spherical suspension construction, (see [Ll1998] and sections  3.4  and   6.3.4  in [Gr2021]  reduces the  case of odd dimensions to that of the  above even ones. QED.

.

{\it \textbf { Remarks.}} \textbf{(i$_C$) } The above mentioned  spherical suspension construction
also allows a reduction of the theorem to the case, where   sphere bundle $S^\bullet\to P$
is {\it trivial}, as follows (compare with he proof of theorem 13.8  in [GL1983]).

Let $S^+\supset S^\bullet$ be an extension of $ S^\bullet \to P$ to a  trivial sphere bundle,
 $S^+=P\times S^{K^+}$, $K^+> K+dim(P)$ and let $\mathbb R_p^{K^++1}=\mathbb R^{K^++1}\supset S^+_p=p\times S^+$ be the Euclidean spaces containing the fibers of this trivial sphere bundle.

Let $R_p^\perp$, $p\in P$,  be the  $(K^+-K^\bullet)$-dimensional family  of  unit vectors $\nu_r\in  T_{-\delta(p)}(\mathbb R^{K^++1})$ normal to  $S^+_p$ 
and pointing  inward  the unit ball in  $\subset \mathbb R_p^{K^++1}$ bounded by $S^+_p$.

Let $S_r\subset S^+_p$, $r\in R_p^\perp$, be   $K$-dimensional subspheres 
$S_r\subset S^+_p$ obtained by intersecting   $S^+_p$ by affine  subspaces of dimension $K^\bullet +1$ in   
$\mathbb R_p^{K^++}\supset S^+_p$, spanned by tangent subspace 
$T_{\delta(p)} (S^\bullet_p)\subset T_{\delta(p)}(\mathbb R_p^{K^++1})\subset \mathbb R_p^{K^++1}$  and vectors $\nu_r$,
 where 
  we allow vectors $\nu_r$ tangent to $S_p^+\subset S^\bullet_p$, where the corresponding $K$-spheres $S_r$ collapse to the point $-\delta(p)\in S^\bullet_p$.  
 
Let  $R^\perp \to P$ be the fibration with the fibers $R_p^\perp$  and let $Q^+\to Q$ be the induced fibration over $Q$.

 The  map $\mathcal F :Q\to S^\bullet$ naturally suspends to a map $\mathcal F^+: Q^+ \to S^{+\ast}$, where $S^{+\ast}\to R^\perp$ is the (trivial!) sphere bundle induced from 
 $S^+\to P$ by the map $R^\perp\to P$,
 and the inequality 3.A for $\mathcal F$ reduces (this is an exercise)  to that for  $\mathcal F^+$.
 
   \textbf{(ii$_C$) } A  similar suspension argument   allows  us to replace $[area<]_{\exists\rtimes}$ by the corresponding inequality with the plane $Sc$   for  the  foliation of  the Riemannian product $Q\times \mathbb R^N$ (regarded as the $\mathbb Z^N$)-covering of 
   $Q\times \mathbb T^N$  in the definition of $Sc^\rtimes$)
    by  the leaves $\tilde X_p\times \mathbb R^N$ endowed with the metrics $g_p^\rtimes$  (lifted from  $Q\times \mathbb T^N$ to  $Q\times \mathbb R^N$). 
    
    This $Q\times \mathbb R^N$ is  mapped to the suspension $S^\bullet\wedge S^N$   via the one point compactification  maps $\mathbb R^N_q\to S^N$  composed with the   $\varepsilon$-scaling map  
    $$Q\times \mathbb R^N \overset {\varepsilon} \to Q\times \mathbb R^N\mbox { for } (q, r)\mapsto (q, \varepsilon r).$$

     (Recall that the spherical suspension   corresponds to adding 
the trivial vector bundle of rank  $N$ to the $(K+1)$-vector bundle associated with $S^\bullet$.)

What is essential is that  

$\bullet$  the  map  $\varepsilon$ preserves the splitting of the tangent bundle 
$T(Q\times \mathbb R^N)=T(Q)\times \mathbb R^N $;

$\bullet$  the  map  $\varepsilon$ is $g^\rtimes$-isometric on the tangent spaces to $Q\times r$;

$\bullet$  the  map  $\varepsilon$ is $\varepsilon$-contracting on all $q\times\mathbb R^N$ with 
respect to the metric $g^\rtimes$.
 
 Therefore,   Llarull’s trace evaluation (4.6 in [Ll1998]) of the curvature term in the twisted Lichnerowitz formula applied to a
  small positive $\varepsilon$ yields  positivity of leaf-wise   Dirac operator twisted with  $\mathbb S_{\mathcal F(p)}^\mp(S_p^\bullet \wedge S^N)$ and the proof follows
  
  Notice that  this argument automatiaclly takes care of odd dimensions and it  needs no  additional twist with almost flat bundles $\mathcal T^\rtimes$.  However, when it comes to  rigidity theorems, such a twist, or rather   
  twisting with families of flat  bundles, makes he proofs easier.

 \textbf{(iii$_C$) }   The spin condition on  $T(\mathcal X)$ can be relaxed to {\it spin of the lift of $T(\mathcal X)$ to the universal covering of $X$}  (as in  \S\S$9\frac{1}{9}$, $ 9\frac {1}{8}$  in [Gr1996] and  section 10 in [Gr2019] and  footnote 78 in  section 2.7 in [Gr2021]) 
but one has {\color{red!55!black}no idea} what to do in the fully non-spin situation.

 \textbf{(iv$_C$) }  The recent development in the Dirac  index theory for manifolds with boundaries (see [CZ2021],  [WXY2021] and references therein) {\color {red!35!black}shows} -- I haven't truly checked this -- that completeness of the fibers $X_p$ can be replaced by a specific (finite!) lower  bound on the distances  from the supports of the differentials of the maps $ \mathcal F_p:X_p\to S^\bullet_p$  to $ \partial X_p$  (and /or to  $ \partial_\infty X_p$).


\vspace {2mm}

Let us now turn to the case $m\leq n$  and let us 
  define  the $\tilde\square^\perp$-multispread of an  $h\in H_{j/k}(Q,Q\setminus \Psi^{-1}(\partial{P}))$  by exhausting $P$ by compact domains $P_i$,
  $$P_1\subset ...\subset P_i \subset ...\subset P,$$
  observing that $h$ can be represented by infinite $j$-cycles $C_i\subset P_i$ for all sufficiently
   large $i$ letting $h(i)=[C_i]\in H_j(\Psi^{-1}(P_i^\circ)$   for the  interiors  $P_i^\circ\subset P_i,$ 
  and setting
  $$\tilde\square^\perp(h)=\limsup_{i\to \infty} \tilde\square^\perp (h(i)).$$


  \textbf {3.D. Area Contraction Conjecture.} {\sf  Let us  allow the fibers $\tilde X_p=\Psi^{-1}(p)\subset Q$ of the submersion $\Psi:Q\to P$
 to have boundaries and  allow  $m\leq n=dim(\tilde X_p)$.  

Let 
  $$h_K\in H_{K/K-m}(Q, \Psi^{-1}(\partial P)),$$ 
and 
   $$ \sigma_N^\rtimes=\sigma_N^\rtimes (q)=Sc(g_p^\rtimes(q))
-\frac{n+N-1}{n+N}\cdot \frac {4(n-m)\pi^2}{\tilde \square^\perp (h_K)^2}.$$ 
 
   Let the map $\mathcal F$ be
   {\it fiberwise area decreasing}  with respect to the metrics 
 $\sigma_N^{\rtimes} \cdot g $ in $Q$ and $\gamma_m$ in $S^\bullet$.

{\it If the fibers $\tilde X_p$ are metrically complete, then 
$$\mathcal F^\ast(\theta^\bullet)(h_K)=0.$$}}
 
(This conjecture points toward  our ultimate goal,  but it may need an adjustment to avoid a possible 
  counterexample in its present formulation, e.g.  obtained  with  a    parametric thin surgery.)

\vspace {1mm}
We prove below a special case of his conjecture by 
 applying 3.B to families of $m$-dimensional submanifolds $Y_p\subset \tilde X_p$  for $m<n=dim(\tilde X_p)$,  where the scalar curvatures of $Y_p$ are bounded from below  according to  the 
$\square_\ast^{\exists\exists}(n,m)$-inequality 2.B.
\vspace {1mm}\vspace {1mm}

\hspace {27mm}{\sc Preparation for Theorem  3.E.}\vspace {1mm}

Let  $P$ and $Q$ be smooth  manifolds, possibly with boundaries,  and let  
$\Psi :  Q\to P$ be   a true fibration, with an oriented, possibly disconnected, $n$-dimensional Riemannian  fiber $\tilde X=(\tilde X, g)$, where the  structure group $\Pi$ of the fibration is countable  and its action
 on $\tilde X$ 
{\it  discrete,  free, orientation preserving} and where, since $Q$ is connected,  the quotient space $\underline X=
\tilde X/\Pi$ is connected.

Here,   the fundamental group $\pi_1(P)$ acts on $\tilde X$ via a  surjective  homomorphism $\pi_1(P)\to \Pi$ and 
$$Q=(\tilde P\times \tilde X)/\pi_1(P)$$
for the Galois action of $\pi_1(P)$ on the  universal covering  $\tilde P\to P$ and the diagonal action of $\pi_1(P)$ on the product $ \tilde X\times \tilde P$.

Thus,   there is a natural   covering map, say  
$$^\sim\hspace {-0,7mm}\Phi: Q\to P\times \underline X= (\tilde P\times \tilde X)/\pi_1(P)\times \pi_1(P).$$

    Given homology classes  $\underline h_m\in H_m  (\underline X,\partial_\infty  \underline X)$,  where   
 $m\leq n=dim(\underline X)=dim(\tilde X)$, and 
 $h_{K-m}\in H_{K-m}(P,\partial P)$, let  
 $$h_K\in H_{K/K-m}(Q, \Psi^{-1}(\partial P))$$
  be the pullback of the  product  
  $\underline h_{K-m}\otimes h_{m}\in H_K(P\times \underline X,
   \partial P\times \underline X\cup P\times \partial_\infty X)$ 
   to $Q$ under the  map $^\sim\hspace {-0,7mm}\Phi: Q\to P\times \underline X.$
  
 Let  $\underline g$ be a  Riemannian metric in $\underline X$ for which $(\underline X,\underline g)$
 is metrically complete (bounded subsets are precompact),  and let 
 $g_p$ be the the  Riemannian metrics in the fibers $\tilde X_p$ induced by the covering maps
 $\tilde X_p=\tilde X\to \underline X$.

  Let $\varphi_i(x)$,  $i=1,...   N$, be    smooth  positive functions on $\underline X$  and
 $$\underline g^\times  =\underline g+\sum_{i=1}^N\varphi(x)^2dt_i^2,$$ 
 be the  metric on $\underline   X\times \mathbb T^N$.

 Recall the     homology class $\underline h_m\in H_m(\underline X, \partial_ \infty \underline X)$
  and let
  $$\underline \sigma_N^\rtimes =\underline \sigma_N^\rtimes(x)= Sc(g^\rtimes)
-  \frac {n+N-1}{n+N}\cdot\frac {4(n-m)\pi^2}{\tilde \square_{\underline g^\rtimes}^\perp (\underline h_m)^2}.$$
 
 \textbf {3.E. Area Contraction  Theorem.} {\sf Let the  above  map 
  $\mathcal F:Q\to S^\bullet$ 
  be  {\it fiberwise area decreasing}  with respect to the metrics 
 $\underline\sigma_N^\rtimes  \cdot g_p  $ in the fibers $\tilde X_p=\tilde  X$ of the fibration 
 $\Psi: Q\to P$
  and with the above normalized spherical metrics  $\gamma_m$ in the fiberes of $S^\bullet\to P$.

 {\it If
  
  $\bullet_{spin} $ \hspace {3mm} either   the manifold  $\tilde  X$ be spin or $m\leq 3$,
   
  $\bullet_{cmplt} $  \hspace {2mm}$\tilde  X=(\tilde  X,g) $ be  metrically complete, 
  
 $\bullet_{\leq 8} $\hspace {6mm} $n=dim (\tilde  X)\leq 8$. 

Then 
$$\mathcal F^\ast(\theta^\bullet)(h_K)=0.$$}}
 
In fact, as we have already  stated, this is a direct corollary of  2.B+3.B. \vspace {1mm}


 \textbf {3.F. Representative Example of 3.D. } {\sf Let $Y$ and $P$ be smooth compact oriented  manifolds  and let  $X$ be a Riemannian manifold homeomorphic to $Y\times Z$, where $Z$ is a compact enlargeable $k$-manifold, e.g. $Z=\mathbb T^k$.
  
{\sf Let $Sc(X)\geq m(m-1) $ for $m=dim(Y).$}

Let $dim(P)=K$,  let  $S^{m+K}\subset \mathbb R^{m+K+1}$ be  the unit sphere and let 
$$f:X\times P\to S^{m+K}$$ 
be a smooth map, which is {\it  area
 decreasing} on all $X=X_p=X\times p\subset X\times P$, $p\in P$.

If the universal covering of $Y$ is  {\it spin} and if $dim(X)\leq 8$, then 
$$f_\ast([Y\times P])=0\in H_{m+K}(S^{m+K})=\mathbb Z,$$}
where $[Y\times P]\in H_{m+K}(X\times P)$ denotes the  homology class of 
the submanifold $Y\times P= Y\times \mathbf 0 \times P\subset X=Y\times \mathbb T^k\times P$,
$\mathbf 0 \in \mathbb T^k$.\vspace {1mm}
ies) 
{\it Remarks.  \textbf{(i$_E$) }  Granted validity of $\square^{\exists\exists}(n>8)$-conjecture 2.C}
 (irrelevance of singularities)   
the proof of 3.D  applies to all $n=m+k$.

But this    proof  breaks down for all dimensions  if we allow mutually non-isometric    fibers  $X_p =(X, g_p)$, 
 (even for 3.E) , since the minimal hypersurfaces and stable 
$\mu$-bubbles  are not continuous in $g_p$.\footnote {We go around    this discontinuity for $m=3$ and $K=1$ in 
the  4d-example 7.C.}   
 
This begs} for a {\it purely Dirac theoretic} proof of 3.D in the spirit of  [CZ2021] and  [WXY2021] that    could yield  more general inequalities, e.g.  formulated in terms of the {\it K-area},\footnote{The only feasible approach to such inequalities lies in the Dirac theoretic  realm but 
  the scalar curvature  geometry of  minimal hypersurface and $\mu$-bubbles depending on parameters is interesting in its own right;   we say a few words about it in section 7.}. 


{\color {red}corrected on May 12 2023}\textbf {3.G. $\bf Rad_{area} (h_\ast/S^K)$-Invariant.} Much (but not all)  of the above can be expressed in the following terms.

 Given   a Riemannian  manifold $X$
and $h_m\in H_m(X, \partial_\infty X)$, let $Rad_{area}(h_m/S^K)$ be  the supremum of the numbers $R>0$  with the following property.

There exits a cellular topological space $P$, a   bundle 
$S^\bullet\to P$ of  $K$-spheres of radii $R$  with a distinguished section $P_\bullet\hookrightarrow S^\bullet$  and a continuous map $F: X  \times   P\to S^\bullet$, such that

$\bullet_{\bullet}$ the map  $F$  sends  
$\partial_\infty  X\to P_\bullet\subset S^\bullet$;  

$\bullet_p$ the  $p$-fiber maps $F_p=F_{|X\times p}: X\to S^\bullet_p=S^K$ of $S$, are smooth  and $C^1$-continuous in $p\in P$;
logycl
$\bullet_{area} $ the maps $F_p$ are area decreasing for all  $p\in P$; 

$\bullet_{top}$ there exists a rational homology class 
$h_{K-m} \in H_{K-m}(P) $ such that  the  $F$-pullback of the Thom class 
$\theta^\bullet\in H^K(S^\bullet, P_\bullet)$    doesn't vanish 
on $h_m\otimes h_{K-m} \in  H_K( X\times P, \partial_X\times P)$,
$$F^\ast(\theta^\bullet)(h_m\otimes h_{K-m} )\neq 0.$$

The proof of 3.C shows that if the universal covering of $X$ is {\it spin}, then 
$$ Sc^{\exists \exists\rtimes} (h_m)\geq\underline  \sigma>0 \implies  Rad_{area} (h_m/S^K)\leq \sqrt \frac {m(m-1)}{\underline\sigma}.$$
This, together with  the   \ $\square^{\exists\exists}(n,m,N)$-inequality  2.B, 
$$Sc^{\exists \exists\rtimes} (h_m) \geq Sc^\rtimes(X))
-  \frac {n+N-1}{n+N}\cdot\frac {4(n-m)\pi^2}{\tilde \square^\perp (\underline h_m)^2},$$
proved for $n=dim(X)\leq 8$
yields  the corresponding  lower  bound on   $Rad_{area} (h_m/S^K)$ in terms of $Sc(X)$ and the $\square^\perp$-spread of $h_m$.

{\it Remark on $\bf Rad_{vol_n} (h_\ast/S^K)$}. The above definition obviously generalizes from $area =vol_2$ to  all $vol_i$, where, for instance, Almgren's max-min argument (details  need checking)  yields the inequality 
$$ Rad_{vol_n} ([X]_\ast/S^K)\leq \sqrt[n]{vol(X)/vol(S^n)}, \mbox { } n=dim(X),$$
and the convex partition argument yields the corresponding $\varepsilon$-waist inequality  under the $\mathbb Z_2$-version of $\bullet_{top}$.
(See [Gu2014] and references therein.)
\section {Injectivity Radii and Maps to Spheres}
{\it Notation and Definitions.}  Let  $P$ be   a   Riemannian  manifold, $\tilde P\to P$ be the universal covering over  $P$  regarded as a principal  fibration with the fiber equal to the 
 fundamental group  $\pi(P)$.

Let 
$$\Psi=\Psi_P:\tilde P^\Delta \to P$$ 
be the fibration over $P$ which has   $\tilde P$  as a fiber,  and which is associated with the universal  covering $\tilde P$ over $P$  for  the deck 
 (Galois) action of the fundamental group  $\pi(P)$ in this fiber.

 Thus, 
 the space    $\tilde P^\Delta$ is equal to  $(\tilde P\times \tilde P)/\pi_1(P)$ for the diagonal   $\pi_1(P)$-action,  where  the two $\tilde P$-factors  of $(\tilde P\times \tilde P)$ play here different roles.  \vspace {1mm}

 Let $\delta :  P \to  \tilde P^\Delta$ be the section of this fibration defined by the diagonal embedding $\tilde P \to (\tilde P\times \tilde P)/\pi_1(P)$. 
 
Let $r=r(p)\geq 0 $ be  a continuous  function on $P$, 

 let $BT(r)=BT(P,r)\subset T(P)$ be the "subbundle"  of $r(P)$-balls 
in the 

tangent spaces $T_p(P)$\footnote{This fails to be a true bundle at the points where $r(p)=0$,}

 let 
 $\tilde B^\Delta (r)\subset  \tilde P^\Delta$  be the 
 "subbundle" of the  $r$-balls $ \tilde B_{\delta(p)}(r(p))$  in the fibers 
 
 $\tilde P_p^\Delta \subset \tilde P^\Delta$, where these  balls degenerate to points  for  r(p)=0. 

Let 
  $r(p)\leq inj.rad_{\delta(p)}(\tilde P_p)$\footnote{If $P$ is  {\it geodesically incomplete}, i.e.  the exponential map at an $p\in P$  is  defined only on certain   open tangent  ball $BT_p(R_p)\subset T_p(P)$, then, by definition, the inequality 
$inj.rad_p\geq r$  signifies  that $r\leq R_p$,  i. e.  the exponential map   $\exp_p:BT_p(r) \to P$ {\it is defined} and it is   diffeomorphism on its image.}  and let $r(p)$  be {\it strictly positive inside} $P$, 
i.e. $r(p)>0$ for  $p\in P\setminus\partial P$,  and observe that these "subbundles" are true ball  bundles over $P\setminus \partial P.$

 In fact,  the inverse exponential maps in the fibers
  $\tilde P^\Delta_{\delta(p)}=\tilde P$ isomorphically send $\tilde B^\Delta_\delta (r)$ to the bundle $BT(r)$,
$$\widetilde {\exp}^{-1}:\tilde B^\Delta_\delta (r)\to BT(r)=BT(P,r), $$
 where the  tangent $r$-balls $BT_p(r(p)) $  are  identified with the tangent balls to the fibers  $\tilde P^\Delta_p $  at $\delta(p) \in \tilde P^\Delta_p$.

Let $S^\bullet=S^\bullet(r)\to P$ be the $K$-sphere "bundle",  $K=dim(P)$ that is  obtained by shrinking the complements  of the open  balls $BT_x(r(x))\subset T_p(P)$ to points and let $P_\bullet \subset S^\bullet$ be the image of the section 
$\delta_\bullet :P\to S^\bullet$ which sends  $p\in P$ to these points.

  If  $p\in \partial P$ and  $r(p)=0$, we agree that the "fiber"  $S_p^\bullet \subset  S^\bullet$  consists of a single point $\bullet_p=\delta_\bullet (p)$.

Let  
 $$\tilde {\mathcal B} : \tilde P^\Delta\to S^\bullet( r)$$
   be    
   {\sf the {\it composition} of  $\widetilde{\exp}^{-1}:\tilde B^\Delta_\delta (r)\to BT(r)\subset T(P)$
  with the   (tautological) quotient map   $T(X)\to S^\bullet( r)$.  }

 Let $E: \tilde P^\Delta\to P$ be the map, where each   fiber $\tilde P_p^\Delta =\tilde P\subset \tilde P^\Delta$  goes to $P$
by the covering map $E_p=\tilde P\to P$, which sends $\delta(p)\mapsto p$.

Let  $X$ be a smooth $m$-dimensional manifold and let $f:X\to P$ be a proper map.

Let $\tilde {Q}_X $ be the product of   $\tilde P^\Delta\underset{E}\to P$ and  $X\underset {f}\to P$
over $P$, i.e.   
$$\tilde {Q}_X =\{ \tilde p^\Delta,x)\in (   \tilde P^\Delta, X\}_{E (\tilde p^\Delta)=f(x)}.$$
(If $X\subset P$, then $\tilde {Q}_X=E^{-1}(X)$.)

Let the map
$$\tilde f^\Delta: \tilde {Q}_X\to \tilde P^ \Delta $$ 
be defined  by the projection $(\tilde p^\Delta,x)\to \tilde p^\Delta\in \tilde P^ \Delta $, thus 
 $E\circ \tilde f^\Delta(\tilde p^\Delta,x)   =f(x)$,  
 and let 
$$\Psi_X= \Psi_P\circ \tilde f^\Delta : \tilde {Q}_X\to  P. $$ 
Observe  that  this   $\Psi_X$ is a fibration over $P$, where the fiber
$\tilde {Q}_{X,p}$ is equal to the covering of $X$ induced
by the map  $f : X\to P$   from the covering $E_p:  \tilde P^ \Delta_p \to P$, and 
that $\tilde {Q}_{X}$  covers the product of $P$ and $X$, where this covering (map), say
$$\Psi_X^X:\tilde Q_X\to P\times X,$$
    is induced from  the covering $\tilde P^\Delta=(\tilde P\times \tilde P)/\pi_1(P)\to P\times P$
by the map $P\times X\to P\times P$  for $(p,x)\mapsto (p,f(x))$.

Let  
$$\tilde {\mathcal F}_{X, r}= \tilde {\mathcal B}\circ\tilde f^\Delta:  \tilde {Q}_X\to S^\bullet (r).$$ 


\vspace {1mm}

\textbf {4.A. Nonvanishing Lemma.}\footnote {Compare with lemma 13.5 in [GL1983].}{\sf Let $P$  and $X$ be
   Riemannian manifolds, possibly with  boundaries and 
$f:X \to P$ be a proper map with the image in the interior  of $P$, i.e.  $f(X)\subset P\setminus \partial P$.
and let $h_m\in H_m(X, \partial_\infty X)$. 

Let $r(x)\leq inj.rad_x (P) $ be a continuous function on $P$, which is strictly positive inside $P$ and 
   which is {\it bounded by distance to $P_0$}   outside a {\it compact subset}    $P_1\subset P$
    $$r(x)\leq dist (x, P_0) \mbox {  for }  x\in P\setminus P_1, $$
where necessarily $P_1\supset P_0.$ 
 
If the $\mathbb  Q$-tensorisation of the image of the class $h_m$ in $P$ doesn't vanish,
$$ 0 \neq f_\ast(h_m)\otimes \mathbb Q\in H_m(P, \partial_\infty P;
\mathbb Q),$$ 
{\it then the pullback  of the    Thom class $\theta^\bullet \in H^k(S^\bullet, P_\bullet;)$ under the map  
$   \tilde {\mathcal F}_{X, r}$ ,
$$ \tilde {\mathcal F}_{X, r}^\ast(\theta^\bullet)\in H^K( \tilde {Q}_X,{\mathcal F}_{X, r}^{-1}(P_\bullet)),$$ 
 doesn't vanish on $H_{K/K-m}(\tilde Q_X,\Psi_X^{-1}(\partial P)).$}}

\vspace {1 mm}

{\it Proof.}  Let $h_{K-m}\in H_{K-m}(P,\partial P)$  for $K=dim(P)$ be a homology class, such that 
the intersection number   $[h_m\cap h_{K-m}]_{ind}\in \mathbb Z$  doesn't vanish (we assume $P$ is connected oriented) 
and let 
$$ \tilde h_K\in H_{K/(K-m)\downarrow \Psi_X}(\tilde Q_X, \Psi_X^{-1}(\partial P\times X\cup P\times \partial_ \infty X))
\subset  H_{K}(\tilde Q_X, \Psi_X^{-1}(\partial P\times X\cup P\times \partial_ \infty X ))$$ 
(for the above   $\tilde Q_X\subset\tilde P^\Delta  \times   X$ and 
$ \Psi_X:\tilde Q_X\to P$)
be  the pullback of 
the product 
$$h_{K-m}\otimes h_m\in H_K (P\times X, \partial P\times X\cup P\times \partial_ \infty X)$$ 
under the natural (covering) map 
$$^\sim\hspace {-0,7mm}\Phi: \tilde Q_X\to P\times X.$$

 The class $\tilde {\mathcal F}_{X, r}^\ast(\theta^\bullet)$ doesn't  change if  we replace $r$ by a smaller continuous positive  function   $r'(x)\leq r(x)$; thus our geometrical
lemma reduces to a {\it purely topological one}, where $r$ is taken  arbitrarily small.

To simplify further,   properly embed $P\subset  \mathbb R^{2K} $ and replace $P$ by  a small regular neighbourhood  $P^+\subset \mathbb R^{2K}$. 

By  the Thom isomorphism,  this reduces the lemma to the case where the $K$-sphere bundle bundle $S^\bullet$ is trivial 
$$S^\bullet=P\times S^K\mbox { (with 2K instead of $K$)}.$$
Here the Thom class $\theta^\bullet$ is induced from the fundamental class of the sphere $S^K$ by the map 
$S^\bullet \to S^K$ and   {\it {\color {blue!72!black} non-vanishing}} of  $\tilde {\mathcal F}_{X, r}^\ast(\theta^\bullet)$
follows by the obvious degrees comparison  argument as 
in [GL1983]. QED.

\vspace {1mm}

\textbf {4.B.} Let  us give geometric conditions that make the class $\tilde {\mathcal F}_{X, r}^\ast(\theta^\bullet)$, hence   homology image  class $f_\ast(X)\in H_m(P)$, {\color {red!42!black} {\it vanish.}}
\vspace {1mm}

{\it Definition of   $Rad^\bullet_p(P)$.} The radius
$$Rad^\bullet_p(P,r)=Rad^\bullet_p(P,r)_{area} \mbox {  for } r\leq inj.rad_p(P)$$
is the supremum of  the numbers $R>0$, such that the open ball $BT_p(r)\subset T_p(P)=\mathbb R^n$  for  $n=dim(P)$  and  $r=inj.rad_p(P)$ 
admits a smooth  $O(n-1)$-{\it equivariant} map  to the $n$-sphere of radius $R$, say
$$\alpha_{x,r}: BT_p(r)\to S^n(R),$$
where $0\in  BT_p(r)$ goes to the south pole and   the boundary sphere $S_p^{n-1}(r)=\partial BT_p(r)$  to the north pole,
and such that the composition of this $\alpha$ with the inverse exponential map that sends  the open  ball $B_p(r)\subset P$ to the $n$-sphere, 
$$\alpha_{x,r} \circ (\exp_p)^{-1}  :B_p(r)\to S^n(R),$$ 
is {\it area decreasing}.


 {\it Example.} Let $a(r, \kappa)$ denote the area of the disk of radius $r$ in the complete simply connected surface with constant curvature $\kappa$, e.g.
$$\mbox {  $a(r, 1)=2\pi(1-\cos r)$,   $a(\pi , 0)=r^2$,  $a(r , -1)=2\pi (\cosh r-1)$,}$$
and let $$R^\bullet(r, \kappa)=\sqrt \frac {a(r, \kappa)}{4\pi} $$ 
be the radius of the 2-sphere with area $a(r, \kappa)$.
 
 Then,  by the textbook  comparison theorem, the inequalities  
 $$\mbox { $sect.curv(P)\leq \kappa$ and $inj.rad_p(P)\geq r$, where    $\kappa \leq \frac {1}{r^2}$,  }$$  
imply that 
 $$rad^\bullet_p(P)\geq R^\bullet(r, \kappa).$$
for all  complete Riemannian $n$-manifold $P$.


\vspace {1mm} 

 \hspace{30mm}{\sc Preparation for Theorem 4.C.} \vspace{1mm}

 Let, as earlier,  $P$ and $X=(X,g)$ be  complete Riemannian manifolds, possibly with  boundaries,  and let $$f:X\to P\setminus \partial P$$
  be a smooth  proper map.
   Let  $U=U_f\subset P$ be the set of points $p\in P$,  which are closer to the the image of $f$,  than to the boundary of $P$,
 $$ U=\{p\in P\}_{dist (p, f(X))\leq dist(p, \partial P)}$$
  and let
   $$\tilde R^\bullet= \tilde R^\bullet_U=  \tilde R^\bullet_f=\inf_{p\in U} {Rad}_{\tilde p}^\bullet(\tilde P ,r=inj.rad_{\tilde p}(\tilde P)).$$
  for the points  $\tilde p$  over $p\in P$ in the universal covering $\tilde P$ of $P$. 
    \vspace{1mm}

Now we are able to state  the {\it {\color {teal} \textbf {main result}}} of the present  paper. \vspace {1mm}
  
  \textbf {4.C.  Injectivity  Radius Theorem.} {\sf Let  the map  $f:X\to P\setminus \partial P$ be   {\it area decreasing} with respect to the metric 
   $Sc_N^\rtimes(X) \cdot g$ in $X$ and the original metric in  $P$ (used  in the definition of $\tilde R^\bullet_f$ and let  
   $$h\in H_m(X,\partial _\infty X).$$

 If the $\mathbb Q$-tensorisation of  the $f$-image  of $h$ in the rational homology group $ H_k(P,\partial_\infty P;\mathbb Q)$   {\it doesn't vanish, }}
   $$f_\ast(h)\otimes\mathbb  Q \neq 0,$$
then,    in the following   $\bullet_1$-$\bullet_3$ cases,
    the $\tilde\square^\perp$-spread of $h$ is related to the above   $\tilde R^\bullet$  by the inequality  
     $$\frac {n+N-1}{n+N }\frac {4(n-m)\pi^2}{(\tilde \square^\perp(h))^2}\leq Sc_N^{\rtimes}  (X)+\frac {m(m-1)}{(\tilde R^\bullet)^2} .\leqno  [\tilde\square^\perp\& \tilde R^\bullet]$$

 {\sl $\bullet_1$   \hspace {1mm}   The universal covering of $X$  is spin and $m=n$,
   
    $\bullet_2$   \hspace {1mm}  the universal covering of $X$  is spin  and $n\leq 8$,
  
  $\bullet_3$ \hspace {1mm} $n\leq 8$ and  $m \leq 3.$}

   \vspace{1mm}

 {\it Proof.} 
  The above example 4.B  shows that he map  $\tilde {\mathcal F}_{X, r}= \tilde {\mathcal B}\circ\tilde f^\Delta:  \tilde {Q}_X\to S^\bullet (r)$  is area decreasing and the proof follows by using   the $\square^{\exists\exists}(n,m,N)$-inequality 2.B 
 (this is where $n\leq 8$ comes from) 
   and    confronting  non-vanishing lemma 4.A with  
 the area contraction theorem 3.C (where the spin condition resides).

  \vspace {1mm}\vspace {1mm}

 \hspace {22mm}{\sc Corollaries, Remarks, Conjectures.}\vspace {1mm}
 

 \textbf {4.D. { \large $\mathbf {(m=n)}$}-Example}. The  inequality  $ [\tilde\square^\perp\& \tilde R^\bullet]$ for $m=n$ translates  as follows.\vspace {1mm}

{\color {blue}{\huge $\star$}\hspace{-0.7mm}$_{m=n}$}  {\sf  Let  $X=(X,g)$ be a  complete oriented Riemannian  $n$-manifold, $n\geq 2$  with {\it spin} universal covering $\tilde X$  and let 
 $P$ be a complete  Riemannian manifold of dimension  $K\geq n$, such that that $sec.curv(P)\leq \kappa\geq 0$ and 
the injectivity radius of the universal covering of $P$  satisfies
 $$inj.rad(\tilde  P)\geq \frac {\pi}{\sqrt  \kappa}.\footnote {{\color {red!55!black} Conjecturally} the inequality  $injrad(\tilde  P)\geq  {\pi}/{\sqrt  \kappa}$ can be relaxed to $conj.rad(\tilde  P)\geq  {\pi}/{\sqrt  \kappa}$ and, possibly,  
  the inequality $sect.curv(P)\leq \kappa$ can be  dropped at the same time. }$$
 
Let  $f:X\to P$ be a smooth quasi-proper  map\footnote{A map $f$ is quasi-proper if 
it extends to a continuous map between the compactified
spaces, from $X^{+ends} \supset X$   to $P^{+ends}\supset P$.}   (e.g. $f$ is proper or it is  locally constant at infinity).
 
Let $Sc^{\exists\rtimes}(X,x)\geq \sigma(x)>0$  for a positive  continuous function $\sigma$  on $X$,
(e.g. $Sc(X)>\sigma(X)$). 

 If     $f$  is {\it area decreasing} with respect to the metric  $\sigma\cdot g$, 
then the rational fundamental homology class of $X$
 {\it is sent by $f$ to zero}  in $H_n(P,\partial _\infty P;\mathbb Q)$,
$$f_\ast[X]_\mathbb Q=0.$$}
  
  (If $\kappa=0$, this   generalizes  theorem 13.8 from  [GL1983] and if $P$ is the unit $n$-sphere,
  this is the {\large $\rtimes$}-stabilized  Llarull-Listing-Zhang theorem.\footnote{This argument, which relies on the trace   inequality 4.6 in [Ll1998], automatically  delivers,  as it was observed in   [Li2010], a sharper result, where the area 
decreasing condition on $f$ is replaced  by the corresponding bound on {\it trace norm} of $\wedge^2 df$  as it is explained in detail  in  section 3.4.1 in  [Gr2021].

Similarly, $m<n$,  one may use  the supremum of the $\wedge^2$-trace  norms   over the restrictions of  the differential $df$   to all $m$-dimensional  subspaces in the tangent bundle $T(X)$.})

    \vspace {1mm}

    \vspace {1mm}
  
  \textbf {4.E.  Generalization and Proof of 1.B}.  {\sf   {\sf Let $X=(X,g)$ be a  Riemannian manifold of dimension $n=m+k$, such that $Sc^{\exists\rtimes}(X)\geq \sigma=\sigma(x)>0$ (e.g. $  Sc(X)\geq \sigma$), 
  let $Z$ be    a    compact connected enlargeable  manifold.   
  
Let 
$$\phi: X \to  Z$$
be a  smooth map,  let be a covering of $X$ induced by $\phi$
 from the universal covering $ \tilde Z\to Z$ and  let $\hat  \phi:\hat X\to\tilde Z$  be the lift of 
 $\phi$ to $\hat X$. 

Let $Y=\hat \phi^{-1}(\tilde z) \subset \hat X$ be the $\hat \phi$-pullback of a generic point $\tilde x\in \tilde Z$,

Let $P$ be a Riemannian   manifold with $sect.curv(P)\leq \kappa\geq 0$ and let $f:\hat X\to P$ be a smooth map, which is strictly area decreasing 
with respect to the metric $${\kappa\over m(m-1)}\cdot  \sigma\cdot g\mbox {  in $X$}.$$.

Let 

 {\it  $\bullet_{spin}$  \hspace {1mm}  either $m\leq 3$ or   the universal covering of  $X$ is {\it spin}, 
 
 $\bullet_{\leq 8}$  \hspace {4mm}    $n=dim(X)\leq 8$.}

\textbf {4.E(i). Compact Case. }  Let $X$ be compact without boundary, let $dim(Z)=k$  
  and let the $f$-image of the fundamental class $h_m=[Y]\in H_m(\hat X)$  doesn't vanish 
  in the rational  homology  $H_m(P;\mathbb Q),$
  $$\hat f_\ast(h_m)\otimes \mathbb Q\neq 0.$$

Then the universal covering of $P$ satisfies
 $$inj.rad(\tilde P)\leq\frac  {\pi}{\sqrt\kappa}.$$

\textbf {4.E(ii). Non-compact Case. }  Let   $dim (Z) =k-1$, 
  let $\hat Z\subset \hat X$  be a closed subset, 
  and let $h_m\in H_m((\hat X\setminus \hat Z)\cap Y) $ be a homology class, such that   the $f$-image of $h_m$ doesn't vanish 
  in the rational  homology  $H_m(P;\mathbb Q).$

 Let $\hat R=dist_{metr} (\hat Z, \partial X)$,  (where, by definition,  $\hat R=\infty$  for geodesically complete  $X$),  let $\hat U=U_r(\hat Z)\subset X$ be the closed (but non-compact )   $r$-neighbourhood of $ \hat Z \subset X$, and let
  $$\hat\sigma =\sigma - {36(n-1) \pi^2\over n\hat R^2}\geq  m(m-1) \kappa.$$
Then the injectivity radius of the  universal covering  $\tilde P$ at the  pullbacks of 
the points $f(x)\in P$ to $\tilde P$ satisfies
$$\inf_{x\in U}   inj.rad_{\tilde f(x)}(\tilde P)\leq \pi \sqrt {1\over \kappa}.$$ }}
 
 {\it Proof.} These propositions are obvious corollaries of 4.C, while the 
 examples   2B(i) and  2B(ii)   for a manifold $X$ follow  from   4.E(i) and 4.E(ii) applied to 
   the  covering map $f: \hat X\to P=X.$

\textbf {4.F. About  Rigidity.} The $\mu$-bubbles   inequalities present  in our argument,  as well as  area inequalities for maps to spheres are accompanied by the corresponding rigidity phenomena  in the extremal cases (see [Ll1998], [Li2010], [Br2011],  [Zh2019], [GZ2021] and section 5.7 in [Gr2021.]). Then bringing these rigidity  properties   together shows that the  inequalities
$$sect.curv (\underline g) \leq \frac {\sigma}{m(m-1)}$$ 
 and 
$$ inj.rad(\tilde (X,\tilde {\underline g}) \geq \pi\sqrt \frac{m(m-1)}{\sigma}$$
together  with  non-vanishing   of    $[Y]\otimes\mathbb Q\in H_m(X,\partial_ \infty X;\mathbb Q)$
imply that the $g=\underline g$ and the  universal covering $\tilde (X,\tilde g)$ is isometric to $S^m(R)\times\mathbb R^k$ for some $R>0$.

\vspace {1mm}


 


\textbf {4.E.}  {\it \textbf {$\tilde {\mathcal  R}_{area} ^\bullet$-Invariant.}} The proof of the  
   injectivity  radius theorem 4.C  remains valid if $R^\bullet$ is replaced 
   by      
  $\tilde {\mathcal  R}_{area} ^\bullet(P)\geq R^\bullet$, which quantifies the 
  $\lambda^2$-enlargeability from [GL1983]  and which is   defined as follows.

Recall the map $\tilde {\mathcal B} : \tilde P^\Delta\to S^\bullet( r)$
and let us agree that the inequality   
$$\tilde {\mathcal  R}_{area} ^\bullet(P)\geq r=r(p)$$
signify that there exists a smooth  map
$$\tilde B: \tilde P^\Delta\to S^\bullet( r)$$
with the following two properties.

$\bullet_{hom}$ the map  $\tilde B$ sends $p$-fibers to $p$-fibers 
$$ \tilde P_p\overset {\tilde B}\to   S_p^\bullet( r(p)), \mbox { } p\in P,$$
such that $\partial    \tilde P_p\cup\partial_\infty P_p$ is sent to the point $\bullet_p\in S^\bullet_p(r(p))$ 
for all $p\in P$, and such that  $ \tilde B$ is homotopic to  $\tilde {\mathcal B} $ by such fiber 
preserving continuous maps $\tilde P^\Delta\to S^\bullet( r)$

$\bullet_{area}$ The map $\tilde B$ is fiber-wise smooth, it is $C^1$-continuous in $p\in P$\footnote {I am not 
certain if this is truly necessary for the validity of 4.C.}  and area decreasing on all fibers  $\tilde P_p$.

{\it Example.} It is shown in (different terms) in [Gl1983] that locally homogeneous spaces $P$ with contractible universal coverings  have  $\tilde {\mathcal  R}_{length} ^\bullet(P)=\infty$,\footnote {"Lenght" instead of "area"  means that the corresponding maps $B: \tilde P_p\overset {\tilde B}\to   S_p^\bullet( r(p))$  are length decreasing.} 
  i.e.  $\tilde {\mathcal  R}_{lengh} ^\bullet(P)\geq r$ for all $r$.

It follows, for instance,  that the Riemannian products of these $P$ with unit spheres satisfy  
$\tilde {\mathcal  R}_{area} ^\bullet (P\times S^l)\geq  \pi-\varepsilon $  for all  $ \varepsilon>0$.


\section {Focal  Radii and Normal  Curvatures}

Let $X =(X,g)$  be a complete   Riemannian $n$-manifold,  possibly with a boundary, e.g. $X=B^{n_1} (1) \times \mathbb R^{n_2}.$

  Let  $Z$ be a
$(k+l)$-dimensional manifold, e.g.  homeomorphic to $\mathbb T^k\times S^l$ or to  $\mathbb R^k\times S^l$ and let $f: Z\to X$ be a smooth immersion.

  Let (compare with the proof   of 1.C)
  $$BT^\perp(R) \subset  T^\perp(Z) $$ 
 be the  $R$-ball subbundle in the normal bundle of $Z$    for  $R< foc.rad(Z)$,  remove arbitrarily small open 
 $\varepsilon$neighbourhood of the zero section $Z\hookrightarrow T^\perp(R)$  from this ball bundle,  let 
 $$P=BT^\perp(R)\setminus U_\varepsilon (Z)$$
 and  endow this $P$ with  the same Riemannian  metric $g$ 
 as $X$, that is induced by the normal exponential map $\exp^\perp :BT^\perp(R) \to X$ from the Riemannian metric $g$ in  $X$.
 
Let $S(r/2)\subset  T=BT^\perp(R) $ be the $(R/2)$-sphere  subbundle,
and let 
 $$X^\perp = X_{r, R-r}^\perp\subset P, \mbox { } r<R/2,$$
   be the bundle of annuli 
    pinched between of the sphere subbundles  
    $$S(r), S(R-r) \subset P.$$

Let  $inj.rad_f(z)  (X)\geq\frac {R}{2}- r$, $z\in Z$, and observe that 
$$dist_P(S(r), S(R-r)) =R-2r$$ in this case.

As in the proof of 1.C, we derive a bound on 
$R=foc.rad(Z)-\varepsilon$, $\varepsilon\to 0$, hence on $foc.rad(Z)$,\footnote {Since we prefer  the {\it closed ball } bundle for $P$ we can't use $R=foc.rad Z$.}
 by  
applying 4.C to the imbedding 
$X^\perp\hookrightarrow P$ as follows.

Let 
$$h_l\subset H_l(Z, \partial_\infty Z;\mathbb Q)$$ 
 be a {\it non-zero  rational homology class}, which has {\it infinite  
$\tilde \square^\perp$-spread} with respect to the induced metric in $Z$.
$$\tilde\square^\perp(h_l)=\infty.$$

Let $m=l+(n-k-l-1)=n-k-1$ and let
$$ h_m\in 
 H_{m} (S(R/2) , \partial_\infty S(R/2)=H_{m} (X^\perp , \partial_\infty X^\perp)= H_{m} (P , \partial_\infty P)$$
 be image  of  $h_l\in H_l( \tilde Z,\partial_\infty \tilde Z;\mathbb Q$ under the ratioanal  {\it Gysin's homomorphism} 
 $$H_l(Z;\partial_\infty Z;\mathbb Q )\to H_m(P;\partial_\infty P;\mathbb Q)$$
  for the $(n-k-l-1)$-sphere bundle $S(R/2)\to Z.$\footnote{The Gysin homomorphism is defined for orientable sphere bundles; if  the bundle $S(R/2)\to Z$   is non-orientable, we pass to an orientable double covering of it.}

Observe that if he   Euler class of the normal bundle $T^\perp(Z)\to Z$, 
$$\chi^\perp_{n-k-l} \in  H_{n-k-l} (Z^\perp),$$ 
vanishes on some
$h_k\in H_k(Z) $, which has {\it  non-zero intersection number with $h_l$},  then {\it $h_m$ doesn't vanish}. This happens, for instance, if 
 $k$ is odd   or if $n-k-l-1$ is even, or if he map $f:Z\to X$ is an {\it embedding} and $H_{k+l}(X)=0$.
 
  Also there often exists a covering $\tilde Z\to Z$, such that the lift of $h_l$ to  $\tilde h_l \in H_l(\tilde Z, \partial_\infty Z;\mathbb Q)$  doesn't vanish, while  the  Euler class of the bundle $\tilde S(R/2) \to \tilde Z$  vanishes on some  $\tilde  h_k $.
  
  For instance if $h_l$ is  the homology  class of the pullback of a generic point under a smooth proper map $Z\to \underline Z$, where  
$\underline Z$ is a compact enlargeable aspherical manifold (e.g. admitting a metric with non-positive curvature)
and $\tilde Z\to Z$ is the covering induced from he universal covering
 of $\underline Z$, then the corresponding $\tilde h_m \in H_m(\tilde P,
\partial_\infty \tilde P;
\mathbb Q)$  {\it doesn't vanish}. 
 
Finally, observe that the $\square^\perp$-spread of the corresponding   
$\tilde h_m\in H_m(\tilde X^\perp,\partial_\infty \tilde X^\perp)=H_m(\tilde S(R/2),\partial_\infty \tilde S(R/2)$ satisfies 
$$\square^\perp(\tilde h_m) \geq dist_P(S(r), S(R-r)) =R-2r$$ 
 for all coverings $ \tilde Z$ of $Z$
and that 
$U\subset P$ defined in 4.C as  
$$ U=\{p\in P\}_{dist (p, X^\perp)\leq dist(p, \partial P)}$$
is equal to the $(\rho, R-\rho)$  annulus bundle for $\rho=r/2$.
 
  \textbf {5.A. Conclusion.} {\sf If the universal covering of $T^\perp(Z) $ is spin and $n=dim(X)=n\leq 8$, }then 
  theorem 4.C applied with this $U$, with 
  $$r=\frac {1}{3}R, \rho =\frac {1}{6} R$$ and with
     $$\tilde R^\bullet\inf_{p\in U} {Rad}_{\tilde p}^\bullet(\tilde P ,r=inj.rad_{\tilde p}(\tilde P)).$$
  for the points  $\tilde p$  over $p\in P$ in the universal covering $\tilde P$ of $P$ shows that 
 $$\frac {n+N-1}{n+N }\frac {4(n-m)\pi^2}{9foc.rad(Z)^2}\geq Sc_N^{\rtimes}  (X)+\frac {m(m-1)}{(\tilde R^\bullet)^2}.\leqno{[\square^\perp \&foc.rad]}$$

{ \sc Discussion.}  
Focal radius is closely related to the {\it maximal  normal curvature} of an immersions $f:Z\to X$,    that is 
  {\it  the supremum of  the $X$-curvatures 
  of the {\it geodesic lines }  in $Z$,} 
i.e. of  the curvatures measured in the Riemannian geometry of
 $X$ of geodesics in 
 $Z$ for the induced Riemannian metric in $Z$.

However,  the  inequality $[\square^\perp \&foc.rad$] yields  a significant  lower  bound on the curvatures of immersions to  "simple"  manifolds $X$, e.g.  to  the unit ball $B^n\subset \mathbb R^n$ (where $foc.rad(Z)=\min(dist(Z, \partial  B^n, 1/curv(Z))$)  only for large $n/m>>8$, where this  inequality remains conjectural.

 In contrast, 
   Anon  Petrunin recently proved he following.
   
  \textbf {5.B.}   {\it \textbf {Asymptotically Optimal  $\sqrt 3$-Inequality.}} 
{\sf If a compact $k$-manifold $Z$   {\it admits no metric with} $Sc>0$,  e.g. $Z$ is enlargeable,  then  immersions $Z\hookrightarrow B^n$
have 
$curv(Z)\geq \sqrt{ 3k/(k+2)}$ 
for all $k$ and $n$. }\vspace {1mm}

(Optimality  of this result is proved  in [Gr2022], where we   construct  immersions 
$Z^k\to  B^{8k^2}$  with   $curv(Z)=  \sqrt{ 3k/(k+2)}+\varepsilon$  and   
 where  we also present  upper and lower bounds on the curvatures of immersions $Z^k\to X^n$ for $n\sim k$.)

\section {Mean Convex  Domains in $\mathbb R^n$}

Let $X$ be an infinite tunnel in the 3d space,  that is a closed subset  $X\subset \mathbb R^3$ diffeomorphic to the cylinder 
$B^2\times \mathbb R$.

\textbf {6.A.  Example}:  {\it \textbf {Fat Mouse in a Narrow Tunnel}.}  {\sf If the mean curvature of  the boundary of $X$ is everywhere $\geq 1=mean.curv(S^2(2))$, 
then a "mouse", which contains an $(1+\varepsilon)$ ball $B^3(1+\varepsilon)$   inside  its body won't be able to crawl trough this tunnel:}
 
 {\it If a connected subset $Z\subset X$ infinitely stretches out to   the both ends of $X$}, {\sf i.e. it is  not contained in
 $B^2\times S$ for a proper  subset $S\subsetneqq \mathbb R$,  {\it then 
the 1-neighbourhood $U_1(Z)\subset \mathbb R^3$ of $Z$  intersects the boundary of $X$.}\vspace {1mm}}

{\it Remarks.}  (a)  The  model  $X$ is the infinite round   cylinder $B^2(1)\times \mathbb R\subset \mathbb R^3$, where $mean.curv(\partial X)=1$, where 
  the optimal  $Z$ is the central line $0\times \mathbb R\subset X$  with $U_1(Z)=X.$

(b) Instructive examples  of  tunnels   with $mean.curv(\partial X)\geq 1$  are  made by    interconnecting infinite chains of  disjoint   balls $B^3_{x_i}(2-\varepsilon)\subset \mathbb  R^3$   of radii $2-\varepsilon$, $\varepsilon >0$, by narrow tubes, such that $X$ contains all these balls,  
  the mean curvature of $X$
 remains $\geq 2-2\varepsilon$ close to the balls and is arbitrarily large everywhere else.\footnote {Similarly to the thin codimension 2  surgery of manifolds with  $Sc>\sigma$   (1.3 in [Gr2021]),    hypersurfaces with $mean.curv>\mu$  admit codimension 1 surgery, (24 in [Gr2017]). 
 
 Namely, let  $V$ be a domain with a smooth boundary in a Riemannian $n$-manifold $U$ and  $Y\subset U$ be a smooth submanifold of dimension $dim(Y)\leq n-2$ with a boundary, 
 such that  $V\cap Y=\partial Y\subset \partial V$, where  
 the intersection between $Y$ and  $\partial V$ is transversal.
 
 Let $\mu(u)$ be a continuous function on $U$, such that $\mu(u)> mean.curv(\partial V) $
 at all boundary points $u\in \partial V$.
 Then the union $V\cup Y$ admits an arbitrarily  small regular neighbourhood 
 $V_+\supset V\cup Y$  with smooth bounary, such that 
 $mean.curv\partial (V_+,u)>\mu(u)$ for all $u\in \partial V_+$.}

\vspace {1mm}

{\it Proof of 6.A.} This   follows from the following proposition together with the
 {\it Gehring-Bombieri-Simon}   linking/filling theorem (see section 8.1 in [Gr1983] and references therein.)

\textbf {6.B. Lemma}.   {\sf Let $X\subset \mathbb R^3$ be a smooth closed connected, possibly infinite,    domain with smooth {\it non-simply connected} boundary. If   $mean.curv (\partial X)\geq 1$,
Then 
the boundary of $X$ contains {\it non-contractible}  closed curves $\Theta\subset \partial X$  of {\it lengths $\leq 2\pi+\varepsilon$} for all $\varepsilon>0$.}

In fact, since $Sc(X)\geq 0$,  theorem 1.1 in [GZ2021] 
applied to large compact parts of $X$, shows that the manifold 
  $X$ contains connected simply connected surfaces $\Sigma_\varepsilon $, for all $\varepsilon>$ with $\partial \Sigma=\Theta\subset \partial X$, which represent 
non-trivial homology classes in $H_2(X, \partial X)=\mathbb Z$ and such that 
$$\int_\Theta mean.curv(\partial X, \theta) d\theta\leq 2\pi+\varepsilon.$$ 
Since $mean.curv(\partial X)\geq 1$, the   lemma follows.

{\it Remarks. } (a) 
The proof in [GZ2021] elaborates on the  following  observation, which applies to all Riemannian 3-manifols $X$ (compare  with [SY1979] and 5.4 in [Gr2014$'$]).

Let  $\Sigma\subset  X$  be a {\it connected stable minimal} surface with boundary  $\Theta \subset \partial X$  and observe that  
the  stability condition implies that   the intrinsic  curvature of $\partial \Sigma \subset \Sigma$ is bounded from below by the mean curvature of $\partial X$ 
at all points  $\theta\in \partial \Sigma\subset X$.
Then,  as in [SY1979],  the stability inequality to the unit  normal field on $\Sigma$ one shows  that
$$\int_\Sigma Sc(X,\sigma)d\sigma+\int_\Theta mean.curv(\partial X, \theta) d\theta\leq 2\pi,$$
which, for  $Sc(X)\geq 0$ and $mean.curv(\partial X)\geq 1$, yields the inequality
$length(\Theta)\leq 2\pi$.

(b) Somewhat paradoxically,   {\it locally length minimizing} ("length-stable") geodesics  $\Theta\subset \partial  X$ may have arbitrarily large lengths and 
the  minimal surfaces $\Sigma_\Theta\subset X$ with boundaries $\Theta$ can't be deformed  to locally minimizing ones with free boundaries by area decreasing homotopies.

In fact, 
  the minimization process of areas of surfaces bounded by long   $\Theta_{min}\subset \partial X$ divides  $\Theta_{min}$  into several shorter curves in $\partial X$.

(c)  If $X$ is a  {\it compact non-simply connected} domain in $\mathbb R^3$  with $mean.curv(\partial X)\geq 1$,  then, by  the rigidity theorem 1.3 in [GZ2021],  the shortest non-contractible
 curve $\Theta_{min}\subset \partial X$ has length$ <2\pi$.

{\color {red!55!black}Probbaly,} if $diam(X)\leq 1$, then $length (\Theta_{min})
\leq 2\pi-0.01$. (Round tori seem good candidates for extremal $X$.)
\vspace {1mm}

Let us  generalize 6.A. to  Riemannian  manifolds $X$,  $Sc(X)\geq \sigma $ of  higher dimensions $n$, where
the boundary of our $X$ is   divided into two parts, 
$\partial X=\partial_{\square}\cup \partial_ {side}$, where these  are smooth domains in $\partial X$, 
which meat along their common boundary, $ \partial  \partial_{\square}=\partial  \partial_ {side}\subset \partial X.$\footnote {Our $\partial {\square}$ corresponds to  $\partial_{eff}$ in [GZ2021].} 

{\it Example.} If  $X$ is the product of a smooth manifold $B$ by the $k$-cube $\square^k=[-1,1]^k$
then $\partial _{side }X=\partial \times \square^k$ and $\partial_{\square}=B\times \partial \square^k$.

(For the above  cylindrical $X=B^2\times \mathbb R$, the side boundary is  $\partial B^2\times \mathbb R$  and 
$\partial_{\square}$   corresponds to the ideal top and bottom of the cylinder.) 

Let $\Psi:X\to \square^k=[-1,1]^k$ be a continuous map which sends $\partial_\square X \to \partial \square^k$  and let,  as in section 2,  $d_i$, $i=1,...,   k$, be the distances between the pullbacks of the opposite faces of the cube $\square^k$ and 
$$d_\square=d_\square(\Psi)=\left (\frac {1}{k}\sum_{i=1}^k\frac {1}{d_i^2}\right)^{-\frac {1}{2}}.$$

Let $Z\subset X$ be a closed subset and $h_k\in H_k(Z, Z\cap \partial_\square X)$ be a homology class 
which doesn't vanish under the homology  homomorphism induced by $\Psi$,
$$0\neq \Psi_\ast(h_k)\in H_k(\square^k, \partial \square^k)=\mathbb Z.$$

\textbf {6.C. $(n,k)$-Mouse  {\color{red!55!black}Conjecture}\&Theorem.} {\sf Let 
$sect.curv(X)\leq \kappa$, $\kappa\geq 0$, 
 let
$$mean.curv(\partial_{side} X)\geq \mu.$$
and let
$$ Sc(X)\geq\frac {n-1}{n}\cdot
  \frac {4k\pi^2}{d^2_\square. } \leqno{\color {blue}[Sc_\square(n,k)]}$$
Then 
$$inj.rad_Z(X)=\inf inj.rad_{z\in Z}(X) \leq r_\kappa,$$
where $r_k$ is equal to  the radius of the ball $B$ in the sphere $S^{n-k}(1/\sqrt \kappa)$, such that  
$mean.curv(\partial B)=\mu$
and where we {\color {blue} prove this if the universal covering $\Tilde X$ is spin and $n\leq 8$.}}

\textbf {6.D. Corollary/Example.} {\sf Let $X\subset \mathbb R^n$ be a closed  domain homeomorphic to 
$B^2\times  Z_0\times \mathbb R$, where $Z_0$ is a compact $(k-1)$-dimensinal manifold, such that    $\tilde  \square^\perp (Z_0)=\infty,$ e.g. $Z_0$ is the $(k-1)$-torus $ \mathbb T^{k-1}$. 

Let $f:  Z_0\times \mathbb R\to X$ be a proper map properly homotopic to the  imbedding
$ Z_0 \times \mathbb R=0\times Z_0\times \mathbb R\hookrightarrow X$.
If $mean.curv(\partial X)\geq n-k-1$, then, provided $n\leq 8$,  
$$dist(f(Z_0\times \mathbb R),\partial X)\leq 1.$$}

In fact, 6.C applies to compact domains $X_1\subset ... \subset X_i\subset .... \subset X$ which exhaust $X$ and where the Euclidean metrics are slightly perturbed to make $Sc(X_i)>\sigma_i$
for $0\leq \sigma_i\to 0$.\footnote{One may also use $g^\rtimes=g_{Eucl}+\phi^2dt^2$ on $X\times \mathbb T^1$,
for 
$\phi(x)=1+\varepsilon\cdot  dist (x, \partial X)$.}  
 \vspace {1mm}

{\it Proof of   6.C.} Our  argument in the proof of the injectivity radius  $\square_\ast^{\exists\exists}(n,m)$-theorem 4.C shows  that 6.C reduces to a  generalization of  the  area contraction  theorem 3. C to families of maps  from    manifolds with mean convex boundaries to balls in  
the spheres $S^N$.

 This generalization reduces to two other propositions 6.E and 6.F below, where 6.E generalizes 2.C to manifolds with mean convex boundaries and which is  resolved for $n\leq 8$, while 6.F generalizes  codimension zero parametric area contraction inequality 3.A, 
 where an   intended proof needs an extension    of  Llarull's (algebraic)  inequality and where
 this inequality for maps from surfaces  follows from [GZ2021].
(This is  sufficient for  the proof of  6.C for  $n-k=2$.\footnote{ 6.C for $n-k=2$ follows from 6.E below, which  itself, as 
in similar cases mentioned 
 earlier,  follows for all $n$  from a   mild generalization of theorem 4.6 from [SY2017].})


 \vspace {2mm}

\hspace {30mm} {\sc Preparations for 6.E.} \vspace {1mm}

  Given  a Riemannian manifold $(Y,g=g_Y)$ with a distinguished domain $\partial_\ast \subset \partial Y$. e.g. $\partial_\ast=\partial Y$, let us incorporate the mean curvature of $\partial$  in the definition of   
  $Sc^{\exists \rtimes}$  from section 2 as follows. 
  
  Let 
  $Sc ( Y\&\partial_\ast)$ denote the pair $(Sc(Y), mean.curv_g(\partial_\ast))$, where $Sc(Y)$ is a  function on $Y$ and $mean.curv_g(\partial_\ast))$ is a function on $\partial_\ast$  and if 
  $$Y^\rtimes=Y\times \mathbb T^N, g^\rtimes =g_Y+ \sum_{i=1}^N\varphi_i^2dt_i^2,$$ then 
    $Sc(Y^\rtimes\&\partial^\rtimes _\ast)$ is the  pair $(Sc(g^\rtimes), mean.curv_{g^\rtimes} ( \partial_\ast\times \mathbb T^N))$.
  
Accordingly,  we introduce the "inequality" 
$$Sc^{\exists\rtimes} (Y\&\partial_\ast)>(\sigma, \mu)$$ for functions $\sigma$ on $Y$ and $\mu$ on 
$\partial_\ast$ as the existence of $g^\rtimes $, such that  $Sc(g^\rtimes)>\sigma $ and  
$mean.curv_{g^\rtimes} ( \partial_\ast\times \mathbb T^N)\geq \mu$,
where this is used below for $\partial_\ast=\partial Y$.
   
   Then, we define $Sc^{\exists\exists\partial\rtimes}$ on   homology classes of Riemannian manifolds $X=(X,g=g_X)$ relative  to a given  $\partial_\ast \subset \partial X$, where we represent such classes 
   $h_m\in H_m (X,\partial_\ast)$ by 
    1-Lipschitz maps from $m$-manifolds $Y$ to $X$, such that $\partial(Y)\to \partial_\ast $,  and where, similarly to section 2,  we use  notation 
   $$Sc^{\exists\exists\partial\rtimes}(h_m)>(\psi, \nu)\mbox  {  for  functions $\psi$ on $X$ and $\nu$ on  $\partial_\ast$}$$
    as the existence of $Y$ and $f$, such that 
    $$Sc^{\exists\rtimes} (Y\&\partial Y)>(\sigma\circ f, \mu\circ f).$$

   Next, let the boundary of $X$ be decomposed as above,   $\partial X=\partial_\square\cup\partial_{ side}$,  
    define $\tilde\square^\perp(h_m)$ for $h_m\in H_m(X, \partial_{side} \cup \partial_\infty X)$
   with  maps $\Psi$ from $X$ to the cube $[-1,1]^k$, $k=n-m$, such that $\partial_\square(X)$ goes to the boundary of the cube and the$\Psi$-pullback of a generic point is homologous to $h_m$.
   
If $n\leq 8$, then  the 
   proof of 2.C extended to $\mu$-bubbles with boundaries in $\partial_{side} X$ (see section in  [Gr2021] and [GZ2021]) 
   yields the following proposition, which remains   {\color{red!55!black} conjectural} for $n\geq 8$.

 \textbf {6.E. $\square_\partial^{\exists\exists}(n,m)$-Theorem}. Let $X$ be a Riemannian $n$-manifold with decomposed  boundary $\partial X=\partial_\square\cup\partial_{ side}$,  let   $X^\rtimes = (X\times \mathbb T^N, g^\rtimes)$ be a 
  $\mathbb T^\rtimes$-extension of $X$ 
   (for $g^\rtimes=g_X^\rtimes =g_X+ \sum_{i=1}^N\phi_i^2dt_i^2$), 
  denote  $\nu^\rtimes  =mean.curv_{g^\rtimes }  (\partial _{side}\times \mathbb T^N) $
and let  $h\in H_m(X, \partial_{side} \cup \partial_\infty X)$, $m=n-k$. 

If $n\leq 8$, then 
$$
Sc^{\exists\exists\partial\rtimes} (h)\geq \left (Sc(X^\rtimes)
-  \frac {N+n-1}{N+n}\cdot  \frac {4k\pi^2}{\tilde \square^\perp (h)^2}, \hspace {1mm}\nu^\rtimes \right), \leqno{\color {blue}[Sc_\partial^{\exists\exists\partial\rtimes}]}$$
 $$
 Sc_{sp}^{\exists\exists\partial\rtimes} (h)\geq \left (Sc(X^\rtimes)
- \frac {N+n-1}{N+n}\cdot  \frac {4k\pi^2}{\tilde \square^\perp (h)^2}, \hspace {1mm}\nu^\rtimes \right), \mbox { if $X$ is spin},
\leqno{\color {blue}[Sc_{\partial,sp}^{\exists\exists\partial\rtimes}]} $$
 $$
 Sc_{ \widetilde {sp}}^{\exists\exists\partial\rtimes} (h)\geq \left (Sc(X^\rtimes)
-  \frac {N+n-1}{N+n}\cdot  \frac {4k\pi^2}{\tilde \square^\perp (h)^2}, \hspace {1mm} \nu^\rtimes \right), \mbox { if $\tilde X$ is spin},
\leqno{\color {blue}[Sc_{\partial,\widetilde {sp}}^{\exists\exists\partial\rtimes}]}$$
 

\vspace {2mm}

\hspace {30mm} {\sc Preparations for 6.F.} \vspace {1mm}

Let $Q=(Q,G)$ be an orientable  $n$-dimensional  Riemannian manifold with a bounary and let $\mathcal X$ be a smooth orientable $m$-dimensional foliation on $Q$ where all leaves  are transversal to the boundary $\partial X$.

Let $\varphi_i(q)$,  $i=1,...   N$, be    smooth positive  functions on $Q$, let 
$$G^\times  =G+\sum_{i=1}^N\varphi_i^2dt_i^2$$ 
and let  
   $Q^\rtimes =Q\rtimes \mathbb T^N=(Q \times \mathbb T^N, G^\rtimes).$

Let $\sigma^\rtimes (q)= Sc(X_q^\rtimes)$, $q\in Q$, where  
$X_q^\rtimes =X_q\rtimes \mathbb T^N\subset Q^\rtimes  $  for  the leaf $X_q\subset Q $ passing through $q\in Q$ and where  $X_q^\rtimes$ is endowed with the Riemannian metric $g_q^\rtimes $ induced from $G^\rtimes$ on $Q^\rtimes\supset X_q^\rtimes.$

Let  $\mu^\rtimes(q)=mean.curv_{g^\rtimes_q}(\partial X_q^\rtimes )$, $q\in\partial Q$.

Let $\underline S^n\left(\sqrt{1/\underline \kappa}\right)$ denote the complete simply connected $n$-space with constant curvature $\underline \kappa$: spherical for $\underline \kappa>0$, Euclidean for $\underline \kappa=0$ and hyperbolic  for $\underline \kappa<0$.\footnote {This notation is  different from that in section 2.1.B$_4$.}  

Let $\underline B= B(\underline\kappa, \underline r)\subset  \underline S^n\left(\sqrt{1/\underline\kappa}\right)$, be the  $\underline r$-ball in this 
 "sphere" and let $\underline \mu=\underline \mu(\underline\kappa,\underline r)$ be the mean curvature of the boundary 
  $\partial \underline B$.
  
  Let $F:Q\to \underline B$ be a smooth map, which sends $\partial Q\to \partial \underline B$  and which is {\it locally constant at infinity}, i.e.  the complement of a compact subset in $Q$  goes to a {\it finite  subset} in $\underline B$. \vspace {1mm}

 \textbf  {6.F. {\color {red!55!black}Conjecture}\&Theorem: Maps from Foliations to Balls.} {\it {\sf Let the manifold $Q=(Q,G)$  be metrically complete and let $\sigma^\rtimes>0$ and $\mu^\rtimes >0$.}
 
If  the differential $dF(q)$ restricted to the leaf $X_q$ and its exterior square  $\wedge^2dF(q)$  satisfy the following inequalities: 
  
 $\bullet_\mu$    \hspace {1mm}   $||dF{|\partial X_q}||<  \frac {\underline \mu }{ \mu^\rtimes (q)}$  for all $q\in \partial Q$,
 
 $ \bullet _\sigma$  \hspace {1mm}   $||\wedge^2dF{|X_q}||\leq \frac {m(m-1)\underline \kappa}{ \sigma^\rtimes(q)}$  for all $q\in Q$;

  \vspace {1mm} 

   \hspace {-6mm}  then the map $F$ has zero degree,  where we {\color {blue} prove this} if  
   
   (i) {\color {blue}$\mathcal X$ 
   is a fibration}    
   
   (ii) the scalar curvature of $\mathcal X^\rtimes$  is  {\color {blue} uniformly positive, $ \sigma^\rtimes \geq \varepsilon >0$,} 
   
    (iii) the  tangent bundle of the lift of the foliation $\mathcal X$  to  the universal covering  $\tilde Q$ of $Q$   is {\color {blue}spin. }}
  
 \vspace{1mm}
  
  {\it About the Proof.}  
If $m=2$ and $N=0$, this follows from the Gauss bonnet theorem  and if $N\geq 1$ a similar proof  seems plausible (but I didn't check this).

But if  $m\geq 3$,  except maybe for $m=3$ and $N=0$, the only available approach to  the proof is via Dirac operators  where either $T(Q)$ or $T(\mathcal X)$  lifted to the universal covering of $Q$ must 
be spin; in fact  if $X$ has no boundary, then    6.F is proved  [SWZ2021], provided either $T(Q)$ or $T(\mathcal X)$ is spin.

For an  individual manifold $X=Q$ with a boundary, where the universal covering $\tilde X$ is  spin, 6.F follows from   {\it Goette-Semllman's theorem} [GS2000]  applied to suitably smoothed double  \DD$X=X\sqcup_{\partial X}X$ (see section 4.3 in  [Gr2021])  mapped to a double of the receiving ball,  
and where  a more satisfactory  proof, which relies on  a relative index theorem and which has an advantage of yielding rigidity of convex domains in symmetric spaces,   is due to John Lott   [Lo2021].

To extend either of these arguments to foliations one needs a  version of 
 inequality 4.6 from [Ll1998], which  allows replacement of the sup-norms in the inequalities $\bullet_\mu$ and $\bullet_\sigma$   by the {\it  normalized   trace norms} 
 
\hspace {-2mm} $||dF|\partial X||\leadsto ||dF|\partial X||_{trace}/(m-1)$ and $||\wedge^dF||\leadsto ||\wedge^dF||_{trace}/m(m-1)$

\hspace {-6mm} (see section 3.4 in [Gr2021], where this is explained for individual manifolds)
and 
  to prove 6.F, one needs only the following. 
 
 \textbf {6.G. Local Comparison  Lemma.} {\sf Let a point $s_\circ\in S^n$  be taken for a {\it  pole} in the sphere a let us call an equatorial subsphere $S^m\subset S^n$  
 {\it radial} if  it is made by geodesic segments radiating  from  the pole $s_\circ\in S^n$, which is equivalent to inclusion   $S^m\ni s_\circ$.

Let $X_\rtimes$ be a Riemannian manifold, let  $x_\rtimes \in X_\rtimes $   and let
$f: X_\rtimes\to S^n$ be a smooth map, such that the differential $df:T(X_\rtimes) \to T(S^n)$  at a point $x_\rtimes \in X_\rtimes $ sends the tangent space $T_{x_\rtimes }(X_\rtimes )$ to the tangent space of a  {\it radial equatorial} 
 (sub)sphere $S^m\subset S^n$,
 $$df(T_{x_\rtimes }(X_\circ)\subset T_{f(x_\rtimes )}(S^m).$$
 
  Let $f': X_\rtimes \to S^m$ be a smooth map,  differential of which at $x_\rtimes$, is equal 
 to that of $f$, i.e. $f'(x_\rtimes )= f(x_\rtimes )\in S^m\subset S^n$ and the linear map
 $d_{x_\rtimes}f': T_{x_\rtimes}(X_\rtimes0\to T_{f(x_\rtimes)}(S^m)$
  is equal to  $d_{x_\rtimes}f$.

  Let    
  $$g_\circ=\phi_1(r)^2dr^2 +\phi_2(r)^2ds^2,\mbox { } r(s) = dist (s, s_\circ),$$
   be a {\it radial}  Riemannian metric on $S^n$,
   which 
 is equivalent in the present case  to {\it invariance of $g_\circ$  under  the orthogonal subgroup $O(n- 1)\subset O(n)$, which fixes the pole $s_\circ \subset S^n$. }

 Let $\mathcal D$ be the Dirac operator on $X$ with the coefficients in the pullback 
 $ f^\ast(\mathbb S^\pm(S^N))$  for   
 the Clifford (spinor)  bundle $\mathbb S^\mp$ over $S^n$,
$$\mathcal D:  C^\infty ((\mathbb S^\mp(X_\rtimes)\otimes  f^\ast(\mathbb S^\mp(S^n)) \circlearrowleft $$
and  
$$\mathcal D': C^\infty ((\mathbb S^\mp(X_\rtimes)\otimes  (f')^\ast(\mathbb S^\mp(S^m)) \circlearrowleft ,$$
be  Dirac operator  with the coefficients in the $f$'-pullback 
 $ (f')^\ast(\mathbb S^\pm(S^m))$,

 Let 
 $$B=\nabla\nabla^\ast :C^\infty ((\mathbb S^\mp(X_\rtimes)\otimes  (f')^\ast(\mathbb S^\mp(S^m)) \circlearrowleft$$
 and 
 $$
 B'=\nabla'(\nabla')^\ast: C^\infty ((\mathbb S^\mp(X_\rtimes)\otimes  (f')^\ast(\mathbb S^\mp(S^m)) \circlearrowleft $$
  be the corresponding (positive) {\it Bochner Laplacians}
  and let us represent the (zero order)  difference operators  $\mathcal D-B$ and 
$\mathcal D'-B'$ by    endomorphisms of the corresponding vector bundles, while keeping the notation 
 $$\mathcal D-B: \mathbb S^\mp(X_\rtimes)\otimes  (f)^\ast(\mathbb S^\mp(S^n)$$
and 
$$\mathcal D'-B': \mathbb S^\mp(X_\rtimes)\otimes  (f')^\ast(\mathbb S^\mp(S^m).$$

  Then {\it the lowest eigenvalue of the (linear selfadjoint operator)  action of $\mathcal D-B$ on the fiber 
   $$[\mathbb S^\mp(X_\rtimes)\otimes  (f)^\ast(\mathbb S^\mp(S^n)]_{x_\rtimes}=\mathbb S^\mp(X_\rtimes)_{x_\rtimes}\otimes \mathbb S_{f(x_\rtimes)}^\mp(S^n)$$
 is equal to  lowest eigenvalue of the operator $\mathcal D'-B'$ on 
  $$[\mathbb S^\mp(X_\rtimes)\otimes  (f)^\ast(\mathbb S^\mp(S^n)]_{x_\rtimes}=\mathbb S^\mp(X_\rtimes)_{x_\rtimes}\otimes \mathbb S_{f(x_\rtimes)}^\mp(S^m)$$}}

 {\it  Proof.} Since  the $g_\circ$-Riemannian (Levi-Cevita) connection  on the normal bundle $T^\perp(S^m)=T(S^n)_{|S^m}\ominus T(S^m)$ is  (obviously) parallel,    
 the $(S^n,g_\circ)$-spin bundle restricted to $S^m$  decomposes into a sum
 of $2^{n-m}$ copies of the  $(S^m, g_\circ)$-spin bundle. 
 
Then  $\mathbb S^\mp(X_\rtimes)_{x_\rtimes}\otimes \mathbb S_{f(x_\rtimes)}^\mp(S^n)$
 decomposes into the corresponding sum of copies of $\mathbb S^\mp(X_\rtimes)_{x_\rtimes}\otimes \mathbb S_{f'(x_\rtimes)}^\mp(S^m)$,
 for $f(x_\rtimes)=f'(x_\rtimes)$ and $d_{x_\rtimes}f=d_{x_\rtimes}f'$ and,
by theorems  II.4.16 and   II.8.17 in [LM 1989] (compare with  section 4 in [Ll1998]),   
 the operator 
 $\mathcal D-B$ on $\mathbb S^\mp(X_\rtimes)_{x_\rtimes}\otimes \mathbb S_{f(x_\rtimes)}^\mp(S^n)$  decomposes accordingly.

 Thus, the two operators have the same eigenvalues  but with different multiplicities.  
 
 \vspace {1mm}
 
\textbf {6.H. Corollary.}  {\sf Let $X_\rtimes $  be a  Riemannian manifold, let
$x_\rtimes \in X_\rtimes $, let    $(S^n,   g_\circ)$ be the $n$-sphere with a radial metric with respect to a given pole in $S^n$  and let 
 $f:X_\rtimes\to \underline S^n$ be as smooth map. 
 
If $sect.curv(g_\circ, x_\rtimes)\geq 0$ and if
   $$||\wedge^2d_{x_\rtimes}f||< \frac {m(m-1)\underline \kappa(x_\rtimes)}{Sc(X_{\rtimes },x_{\rtimes })  }$$
 then, in the cases  $\bullet_{cnst}$  and $\bullet_\circ$ below,

  the operator 
 $$(\mathcal D-B)_\rtimes: [\mathbb S^\mp(X_\rtimes)\otimes  (f)^\ast(\mathbb S^\mp(S^n)]_{x_\rtimes} \circlearrowleft$$
is positive, i.e.  the lowest eigenvalue of $(\mathcal D-B)_\rtimes $  is >0.

 $\bullet_{cnst}$  The sectional curvature of   the metric $g_\circ$ at the point  $f(x_\rtimes)\in S^n$ is {\it   constant} on the set of  2-planes in $T_{f(x_\rtimes)}(S^n)$.   

$\bullet_\circ$
  The image $d_{x_\rtimes}f ( T_{x_\rtimes })\subset T_{f(x_\rtimes)}(S^n)$
 is {\it contained in the tangent subbundle} of a radial sphere $S^m\subset S^n$.}

\vspace{1mm}

{\it Proof.} The case  $\bullet_{cnst}$  follows from the   trace  inequality
4.6 from [Ll1998], (also see proposition 2 in [Li2010]).

The above lemma reduces the  general  case  of  $\bullet_\circ$, where $n\geq m$, 
to $n=m$, where inequality   1.11 in [GS2000], applies. (Also see proposition 1 in  [Li2010].)

{\it Remark. }  In both cases $\bullet_{cnst}$ and  $\bullet_\circ$, the $f$-pullback of the curvature operator of $g_\circ$ to $X_\rtimes$ is diagonalizable at $x_\rtimes$ and 6.G follows by the argument from  [Ll1998].

\vspace {1mm}

{\it Conclusion of the  Proof of 6.F.}
The doubling\&smoothing argument in   section 3.5 of [Gr2021]\footnote {See [GL1980], [Mi2002]
[BMN2011], [Gr,2014$'$] [BH2021] for various versions of this argument.}
reduces 6.F to the corresponding 
property of maps from foliations on manifolds {\it without boundaries} to spheres with 
radial metrics, where, in the spin case,  the argument  in the proof of the mapping theorem 1.2  from  [SWZ2021] applies.

However,  we need an extension of  [SWZ2021]-arguments to a
  $ \mathbb T^\rtimes $-stabilized scalar curvature  $ Sc^\rtimes (\mathcal X)=Sc(\mathcal X\rtimes \mathbb T^N)$ , instead of the   plain
 $ Sc(\mathcal X)$ and also to where  the universal covering  $\tilde {\mathcal X}$, rather than  $\mathcal X$ is spin.

For this reason,  {\color {blue} we     claim 6.F only for fibration with uniformly positive $Sc^\rtimes$,  but this is sufficient for the  
  proof of 6.C for $n\leq 8$.}

 \vspace {2mm}

 \section {A Few Words on Foliations}

 Let $Q$  be a compact  smooth orientable manifold  and      $\mathcal X$ a smooth $n$-dimensional  Riemannian foliation\footnote {{\it "Riemannian" } refers to a  smooth positive quadratic form on the tangent bundle $T(\mathcal X)\to Q$.} $Q$ with leaves   $ X=X_q\subset   Q$, the scalar curvatures of which are bounded from below by $ Sc(X_q) >0$, which is also written as $Sc(\mathcal X)>\sigma$. 
 
 Let $h\in H_K(Q)$ be a homology class and  let $A\to S^K$ be a smooth  leaf-wise area decreasing map,  which doesn't vanish at $h$  under the  homology  homomorphism
 induced by $A$, 
  $$A_\ast(h)\neq   0\in H_K(S^K)=\mathbb Z.$$
 
 \textbf {7.A.{ \color {red!55!black}Problem}.} {\sf Given $\sigma_-\leq  m(m-1)$ for $m=n-k$  and $k=dim(Q)- K$,  find topological and/or geometric conditions on $h$
that would imply an inequality 
 $$\sigma\leq \sigma_-.  $$}

{\it Remark.} One could derive some  (non-sharp) inequalities of this kind,   from the corresponding 
(mainly conjectural) ones for  individual manifolds with $Sc^\rtimes \geq \sigma$, if one could proof 
the following.

 \textbf {7.B.{ \color {red!55!black} Conjecture}.} {\sf The inequality  $Sc(\mathcal X)\geq \sigma,$ and even $Sc^{\exists \rtimes }(\mathcal X)\geq \sigma,$  imply that $Sc^{\exists\exists\rtimes}(Q)\geq \sigma$.}

(This is proven in [Gr2021$'$] for codimension 1 foliations $\mathcal X$ on manifolds $Q$ with $dim(Q)\leq7$.)

  \textbf {7.C. 4d-Example.} {\sf Let $Q$ be a compact  orientable 4-manifold and $\mathcal X$
   an orientable 
  3-dimensinal Riemannian  foliations with $Sc(\mathcal X)>2$ and such that }

  $\bullet_{<\infty}$ the leafs $X\subset Q$ of $\mathcal X$ have at most finitely many ends,\footnote{It is {\color {red!55!black}  
 unclear} if  this condition is needed.}
   
   $\bullet_{\neq 0}$ there exists a leaf $X_0\subset Q$ of $\mathcal X$, and a homology class  
   $h\in H_2(X_0)$ such that the image of $h$ in $H_2(Q)$, the image of which in  $H_2(Q)$
   under the  homology inclusion homomorphism $H_2(X_0)\to H_2(Q)$ is {\it is non-torsion.}
 
 Then there exists a {\it non-torsion homology  class} $ h_\circ\in H_3(Q)$, such that all smooth 
 {\it leaf wise area decreasing} maps  $f:Q\to S^3$   send $h_3$ to {\it zero}.
 $$f_\ast(h_3)=0\in H_3(S^3)=\mathbb Z.$$ 
 
 Moreover, this  conclusion holds if  the inequality $Sc(\mathcal X)>2$  is {\it relaxed to
$Sc^\rtimes (\mathcal X)>2$.}

 {\it Proof.} Start by recalling the following.
  
    \textbf {7.D. Lemma.} {\sf Let $X$ be a complete orientable Riemannian 3-manifold with 
    $Sc(X)> \sigma+2\varepsilon >0$, possibly, with a  uniformly mean convex boundary (i.e.  $mean.curv(\partial X)\geq \mu>0$),  where $\sigma, \varepsilon>0$,  and let $Y_0\subset X$ be  a compact connected oriented  surface with 
    $area (Y_0)\geq \frac {8\pi}{\sigma}$.
  
  Then there exists a one parameter deformation      $Y_t$ of $Y_0$, $0\leq t\leq 1$, where all  $Y_t$ are piecewise smooth 2-cycles, which {\it continuously depend on $t$ with respect to the flat topology} and such that 
    
    (i) $area (Y_t)\leq  area(Y_0)$ for all $t$; 
    
    (ii) $Y_1\subset X$ is a smooth  surface  with   $area (Y_1)\leq \frac {8\pi}{\sigma+\varepsilon}$.}
     
   In fact, $Y_0$ admits an area decreasing deformation to a stable  minimal surface $Y_{min}$, where  
the the second variation formula  and the  Gauss-Bonnet theorem, 
  imply that 
$$\int_{Y_{min}} Sc(X,y) dy  \leq \int_{Y_{min}} Sc(Y_{min},y)dy \leq 8\pi. \footnote {The  second variation formula was brought to the needed   form in  [SY1979], see section  2.5  in    [Gr2021]  for  details.}\leqno{[SYGB]}$$
  
  Furthermore, the same argument applies to  $Sc^\rtimes (X)> \sigma+2\varepsilon >0$ instead of  $Sc(X)> \sigma+2\varepsilon >0$ with a use of the   Zhu  area  inequality  (see [Zh2019] and section 2.8 in [Gr2021]).
 
 {\it Sketch of the Proof of 7.C.}  Let $\nu$ be a  nonvanishing vector field on $Q$ nowhere tangent  to the 
 leaves of $\mathcal X$.
 
 Let $Y_0\subset X_0$ be a  compact smooth connected surface,  with {\it non-torsion} fundamental  homology  class $[Y]\in H_2(Q)$, where, by the Lemma, one  may be assume that  w
 $area(Y_0)\leq 4\pi-2\varepsilon$. 
 
 Continuously move $Y_0$  in the direction of  $\nu$ by small distance $\delta>0$ locally away from $X_0$ to another leave  $X_1\subset Q$  locally "downstream" of $X_0$ relative to $\nu$  while keeping $area \leq 4\pi-\varepsilon$ all along.
 Then   move the resulting surface within $X_1$ to get $Y_1\subset X_1$  with $area (Y_1)\leq 2-2\varepsilon$.  
 
 Keep doing this  and   simultaneously regularize  the surfaces, thus obtaining a sequence    of surfaces $Y_2,..., Y_i,...$ in different  leaves, where the diameters of $Y_i$ with respect to the induced Riemannian metrics and the curvatures of $Y_i$ in the leaves are uniformly bounded.
  
 It follows that there are arbitrarily large $i$ and $j$,  such that  $Y_i$ and $Y_{i+j}$ come arbitrarily close together in $Q$. 
 
Join the two by a narrow cylinder, and  obtain a 3-cycle $\mathcal Y_{ij}$ sliced by surfaces (some of which may be singular)  of areas  
  $\leq 4\pi-\varepsilon$.
  
Now,  the condition  $\bullet_{<\infty}$ and the area bounds on (almost minimal)  surfaces  ($\mu-bubbles$) for  3-manifolds  $X$ with   $Sc^\rtimes >\sigma>0$, shows that the leaves  $X$ can be exhausted by compact domains  
$$B(1)\subset...  \subset B(2)\subset... \subset X$$
with 
$$area(\partial k)\leq const.$$
 Hence, these $B_k$ can be completed to 3-cycles $C_k\supset B_i$ with 
 $$vol_3(C_k\setminus B_k)\leq const'.$$

It is also clear that some leaf $X$  transversally  intersect $C$ over a  very large (of order of $j$) number $l$ of our surfaces $Y_{i_1},...,Y_{i_l}$, which implies that the homology class of 
intersection $\mathcal Y_{ij}\cup C_k$ for large $j $  and $k$ represent a {\it non-zero} multiple 
of $h\in H_2(Q)$. 
Hence, the class $[\mathcal Y_{ij}]\in H_3(Q)$ is non-torsion as well.

Finally,  since $\mathcal Y_{ij}$ is sliced by surfaces  with areas $<\4\pi$, area decreasing maps $\mathcal Y_{ij}\to S^3$ have {\it degrees zero}  by Almgren's  waist inequality (See [Gr2003], [Gu2014]) and the proof is achieved by taking $h_\circ= [\mathcal Y_{ij}]$ with the above (large) $i$ and $j$.

{ \it Sub-Example.}  Let $\underline {\mathcal X}$ be a usual geodesic foliation on the 2-torus $\mathbb T^2$ and $\mathcal X_0$ be the corresponding foliation on 
$Q= \mathbb T^2\times  S^2$ with the leaves $X=\underline X\times S^2$, where
these $\underline X\subset \mathbb T^2$ may be circles or lines.

Observe that the  product  metric $G_0=dt^2+ds^2$ in this foliation has $Sc=2$, and that all non-zero 
$h\in H_3(Q)=\mathbb Z\oplus \mathbb Z$ 
allow  area non-increasing maps $f=f_h:Q\to S^3$ with $f_\ast(h)\neq 0$.

Let $\underline Y\subset \mathbb T^2$ be a closed geodesic transversal to $\underline X$ and 
let $Y=\underline Y\times S^2\subset Q$ and let 
$$G_\varepsilon=dt^2+(1+\varepsilon(t))ds^2,$$
where $\varepsilon (t)$, $t\in \mathbb T^2$, is a smooth function which is  negative away from a small $\delta$-neighbourhood $U_\delta(\underline Y)\subset \mathbb T^2$ and is small positive in this neighbourhood.

Apparently,  (I haven't check this carefully) there are functions  $\varepsilon(t)$ of this kind  arbitrarily  
$C^\infty$-close to zero, such that 
$Sc(\mathcal X, G_\varepsilon)>2$ and where 
all leaf-wise $G_\varepsilon$-area decreasing maps $f:Q\to S^3$ send the class $[Y]\in H_3(Q)$  to zero.

Yet all $ h'\in H_3(Q)$,  which are {\it not multiples of} $h$, do allow  leaf-wise $G_\varepsilon$-area decreasing maps $f':Q\to S^3$ with 
$f'_\ast(h')\neq 0.$\vspace{1mm}

  \textbf {7.C. {\color {red!55!black} Question}}. {\sf Is there  a meaningful version of lemma 7.D for higher dimensions  and codimensions? }
  
  One knows, for instance,  that 
  
  {\sf if a compact Riemannian manifold $X$ with $Sc(X)\geq 2$ is homeomorphic to $S^2\times \mathbb T^k$, 
  and if $k\leq 6$, 
  then 
     the homology class of the 2-sphere $S^2=S^2\times \{t\}\in X$ is representable by a smooth surface with
     $area\leq 4\pi$} (see [Zh2019]). 
     
     But  if $n=dim(X)\geq 5$ (I am {\color {red!55!black} uncertain} about  $n=4$) a  manifold $X$ with $Sc(X)\geq 2$ may contain locally minimal 2-spheres    homologous to $S^2$ of {\it arbitrary  large} areas.\footnote{Such an $X$  can be obtained by  {\it thin surgery} (section 1.3  in  [Gr2021]) applied to the product $S^2\times \mathbb T^k$, $k\geq 3$,  in a small neighbourhood of $S^2\times \mathbf 0\subset S^2\times \mathbb T^k$.}  
     
    What is more promising   in this regard is the geometry of the functional 
    $$\mbox {$Y\mapsto \frac {1}{Sc^\rtimes  (Y)}$ 
     on the space $\mathcal Y=\mathcal Y(X)$,}$$
       of closed cooriented  smooth (singular?)  {\it hypersurfaces} $Y\subset X$ endowed  with induced Riemannian metrics.    
 
 {\color {red!55!black} Hopefully,} the lower bound on ${Sc^\rtimes} (X)$,   does have non-trivial influence on the {\it topologies}\footnote{Defining  an adequate topology in $\mathcal Y$ is  quite an issue here.} of the {\it sublevels}   $\mathcal Y_{\leq a}\subset \mathcal Y$,  
  of this (or a similar) functional,  where  $\mathcal  Y_{\leq a}$ is the space  
 of  hypersurfaces  $Y\subset X$   
  with $Sc^\rtimes(Y) \geq a^{-1}$, 
  and where  
  we are concerned with {\it connectedness} of   $\mathcal  Y_{\leq a}$  for  $a=Sc^\rtimes (X)$, as well as  with the {\it  images of the inclusion homology homomorphisms} 
  $H_\ast(
  \mathcal Y_{\leq a}) \to  H_\ast(\mathcal Y_{\leq a+b}) $  where $a$ and $b/a$  tend to infinity.


 \section {References} ${}$

[Ba1996] {\sl Disques extrémaux et surfaces modulaire.}
Annales de la Faculté des sciences de Toulouse : Mathématiques, Serie 6, Volume 5 (1996) no. 2, pp. 191-202.

\vspace {2mm}

[BH2009] M. Brunnbauer, B. Hanke, {\sl Large and small group homology}, J.Topology 3 (2010) 463-486. \vspace {2mm}

[BH2021]  C. B\"ar, B. Hanke, {\sl Boundary conditions
for scalar curvature}

 arXiv:2012.0912
\vspace {2mm}

[BMN2011] S. Brendle, F. Marques, and A. Neves, {\sl Deformations
of the hemisphere that increase scalar curvature}, Invent. Math. 185,
175-197 (2011)
\vspace {2mm}

[Br2011] S. Brendle, {\sl Rigidity phenomena involving scalar curvature,}

arXiv:1008.3097\vspace {2mm}

[Ce2020] S. Cecchini, {\sl A long neck principle for
Riemannian spin manifolds with positive scalar curvature,} 

arXiv:2002.07131v1

 \vspace {2mm}
 
 [CL2020]  O. Chodosh, C. Li, {\sl Generalized soap bubbles and the topology of manifolds with positive scalar curvature},  arXiv:2008.11888v3\vspace{1mm}  \vspace{2mm}

[CZ2021]S. Cecchini, R Zeidler,     {\sl  Scalar curvature and generalized Callias operators}

\url{https://uni-muenster.sciebo.de/s/QkocJItJAdY2YJ6}

 \vspace {2mm}
 
 [F-CS1980] D. Fischer-Colbrie, R. Schoen, {\sl The structure
of complete stable minimal surfaces in 3-manifolds of non-negative scalar
curvature}, Comm. Pure Appl. Math., 33 (1980) 199-211. \vspace {2mm}

  [G1963] L. Green,  {\sl  Auf Wiedersehensflachen}
 Annals of Mathematics Vol. 78, No. 2 (1963), pp.
289-299
  \vspace {2mm}

  [GL1980]  M.Gromov, B Lawson, {\sl Spin and Scalar Curvature in the
Presence of a Fundamental Group} Annals of Mathematics, 
c111 (1980), 209-230.
\vspace{2mm}

  [GL1983] M.Gromov, B Lawson, {\sl Positive scalar curvature and
the Dirac operator on complete Riemannian manifolds}, Inst. Hautes Etudes
Sci. Publ. Math.58 (1983), 83-196.

 \vspace {2mm}

  [Gr1983] M.Gromov, {\sl    Filling Riemannian manifolds,} J. Differential Geom. 18
(1983), no. 1, 1-147.

   \vspace {2mm}

   [Gr2003] M. Gromov, {\sl Isoperimetry of waists and concentration of maps}. Geom. Funct.
Anal. (GAFA), Vol. 13, No. 1, (2003), 178-215.

\vspace {2mm}

 [Gr2014]  M. Gromov, {\sl Manifolds: Where do we come from? What are we? Where are we going? }  In: The Poincaré Clay Math. Proc., 19, Amer. Math. Soc., Providence, RI, 2014.\vspace {2mm}
    
 [Gr2014$'$]  M. Gromov,  {\sl Dirac and Plateau billiards in domains with corners}
Open Mathematics 12.8 (2014): 1109-1156. \url {http://eudml.org/doc/269152}

\vspace {2mm}

[Gr2017]  M. Gromov, {\sl 101 Questions, Problems and Conjectures around Scalar 
Curvature}

\url{https://www.ihes.fr/~gromov/category/positivescalarcurvature/}

\vspace {2mm}

[Gr2018] M. Gromov, {\sl Metric Inequalities with Scalar Curvature}, Geometric
and Functional Analysis Volume 28, Issue 3, pp 645-726.

\vspace {2mm}

[Gr2019] M.Gromov, {\sl Mean Curvature in the Light of Scalar Curvature,}

arXiv:1812.09731v2  \vspace {2mm}

[Gr2020] M.Gromov, {\sl No metrics with Positive Scalar Curvatures
on Aspherical 5-Manifolds, }  ¯arXiv:2009.05332  \vspace {2mm}

  [Gr2021] M.Gromov, {\sl Four Lectures on Scalar Curvature}  	
arXiv.1908.10612

   \vspace {2mm}
   
   [Gr2021$'$] M.Gromov, {\sl Torsion Obstructions to Positive Scalar Curvature,}
   
   	arXiv:2112.04825

 \vspace {2mm}
  
   [Gr2023] M.Gromov, {\sl In preparation.}
   
    \vspace {2mm}
 
  [GS2000] S. Goette and U. Semmelmann, {\sl Scalar
curvature estimates for compact symmetric spaces}. Differential Geom. Appl.
16(1):65-78, 2002.

 \vspace {2mm}

 [Gu2014]  L. Guth,  {\sl The waist inequality in Gromov’s work}. In The Abel Prize 2008-2012.
Edited by Helge Holden and Ragni Piene. Springer, (2014), 181–195.
    
\vspace {2mm}

 [GZ2021]  M.Gromov, J. Zhu, {\sl Area and Gauss-Bonnet inequalities with scalar curvature,}
 arXiv:2112.07245

 \vspace {2mm}
 
     [LB2021]  C.  LeBrun, {\sl On the Scalar Curvature of 4-Manifolds,}

arXiv:2105.10785 \vspace {2mm}

     [Li2010])   M. Listing, {\sl Scalar curvature on compact
symmetric spaces}. 

arXiv:1007.1832, 2010.
    
   \vspace {2mm}
   
    [Ll1998] M. Llarull, {\sl  Sharp estimates and the Dirac
operator}, Mathematische Annalen January 1998, Volume 310, Issue 1, pp 55-
71.

\vspace {2mm}

[LM2021] Y. Liokumovich, D. Maximo,  {\sl Waist inequality
for 3-manifolds with positive scalar curvature.}
 	
	arXiv:2012.12478 

\vspace {2mm}

[Loh2018] J. Lohkamp, {\sl Minimal Smoothings of Area Minimizing
Cones}, arxiv.org/abs/1810.03157\vspace {2mm}
  
[Loh2018$'$]  J. Lohkamp, {\sl Contracting Maps and Scalar Curvature},

arXiv:1812.11839
   \vspace {2mm}

    [Lo2020] J. Lott, {\sl Index theory for scalar curvature on manifolds
with boundary}, arXiv:2009.07256    \vspace {2mm}

    [Mi1974]   A.Mishchenko,  {\sl Infinite-dimensional
representations of discrete groups, and higher signatures,} Izv. Akad. Nauk
SSSR Ser. Mat., 38:1 (1974), 81-106; Math. USSR-Izv., 8:1 (1974), 85-111 c \vspace {2mm}

 [Mi2002] P. Miao, {\sl Positive mass theorem on manifolds admitting
corners along a hypersurface, }  Adv. Theor. Math. Phys., 6 (6): 1163–1182, 2002.    \vspace {2mm}
    
  [MN2011] F. Marques, A. Neves,{\sf  Rigidity of min-max minimal spheres in three manifolds}, 
  
arxiv.org/pdf/1105.4632.

   \vspace {2mm}
  [Sm2003]  N. Smale, {\sl  Generic regularity of homologically
area minimizing hyper surfaces in eight-dimensional manifolds}, Comm.
Anal. Geom. 1, no. 2 (1993), 217-228   \vspace {2mm}
        
      [Ri2020]    T. Richard {\sl On the 2-systole of stretched enough
positive scalar curvature metrics on $S^2 \times  S^2$}, arXiv:2007.02705v2.

           \vspace {2mm}
         
         [SWZ2021] G. Su, X.Wang, W. Zhang,
{\sl Nonnegative scalar curvature and area decreasing maps on complete foliated
manifolds}. arXiv:2104.03472
       \vspace {2mm}
   
    [SY1979].   R.Schoen, S.T Yau, {\sl  Existence of incompressible
minimal surfaces and the topology of three dimensional manifolds of non-negative
scalar curvature,} Ann. of Math. 110 (1979), 127-142
      \vspace {2mm}

[SY1979$'$] R. Schoen and S. T. Yau, {\sl On the structure of manifolds
with positive scalar curvature}, Manuscripta Math. 28 (1979), 159-183.   \vspace {2mm}

[SY2017] R. Schoen and S. T. Yau, {\sl Positive Scalar Curvature
and Minimal Hypersurface Singularities}. arXiv:1704.05490 ¯
   \vspace {2mm}

  [WXY2021] J. Wang, Z. Xie, G. Yu, {\sl An index theoretic proof of Gromov’s
cube inequality on scalar curvature}, ¯arXiv:2105.12054    \vspace {2mm}

   [Ze2019] R. Zeidler, {\sl Band width estimates via the Dirac operator}, 
   
   arXiv:1905.08520v2   \vspace {2mm}

[Ze2020] R. Zeidler, {\sl Width, largeness and index theory}. 

arXiv:2008.13754   \vspace {2mm}
    
    [Zha 2020] W. Zhang, {\sl  Nonnegative Scalar Curvature and
Area Decreasing Maps}, arXiv:1912.03649v3.   \vspace {2mm}

[Zha(deformed Dirac) 2021] W.Zhang, {\sl Deformed Dirac Operators and the
Scalar curvature}, to appear in Perspectives in Scalar Curvature, volume edited by M. Gromov and B. Lawson, World Scientific, 2022.      \vspace {2mm}

    [Zhu2019] J. Zhu,{\sl  Rigidity of Area-Minimizing 2-Spheres in 
n-Manifolds with Positive Scalar Curvature}, arXiv:1903.0578 \vspace {2mm}
    
    [Zhu2020] J. Zhu, {\sl Rigidity results for complete manifolds with
nonnegative scalar curvature}
arXiv:2008.07028.

 
  \end {document}